\documentclass[psamsfonts]{amsart}

\usepackage{amssymb,amsfonts,amsthm}
\usepackage[all,arc]{xy}
\usepackage{enumerate}
\usepackage{mathrsfs}
\usepackage{hyperref}
\usepackage{dsfont}
\hypersetup{
pdftitle={Paper Xp for Chaos},
pdfsubject={Mathematics, Harmonic Analysis, Functional Analysis},
pdfauthor={},
pdfkeywords={}
}
\usepackage{bm}
\usepackage{pgfplots}

\allowdisplaybreaks

\newtheorem{theorem}{Theorem}[section]
\newtheorem{corollary}[theorem]{Corollary}

\newtheorem{remark}[theorem]{Remark}
\newtheorem*{remark*}{Remark}

\numberwithin{equation}{section}
\numberwithin{theorem}{section}

\theoremstyle{theorem}
\newtheorem{ltheorem}{Theorem}

\newcommand{\R}{\mathbb{R}}
\newcommand{\C}{\mathbb{C}}
\newcommand{\Z}{\mathbb{Z}}

\newcommand{\F}{\mathbb{F}}

\newcommand{\suma}[1]{\sum\limits_{#1}}
\newcommand{\Suma}[2]{\sum\limits_{#1}^{#2}}

\newcommand{\cexp}{\mathsf{E}}
\newcommand{\laplacian}{\Delta}
\newcommand{\gradient}{\mathrm{D}}

\newcommand{\G}{\mathrm{G}}
\newcommand{\V}{\mathrm{\mathcal{L}(G)}}
\newcommand{\Hcal}{\mathcal{H}}
\newcommand{\A}{\mathcal{A}}
\newcommand{\B}{\mathrm{B}}
\newcommand{\Ss}{\mathsf{S}}

\newcommand{\dem}{\noindent {\bf Proof. }}
\newcommand{\fin}{\hspace*{\fill} $\square$ \vskip0.2cm}

\begin{document}

\null

\vskip-5pt

\null

\begin{center}
{\huge Trigonometric chaos and $\mathrm{X}_p$ inequalities I} \\ {\Large ---Balanced Fourier truncations over discrete groups---}

\vskip15pt

{\sc {Antonio Ismael Cano-M\'armol, \\ Jos\'e M. Conde-Alonso and Javier Parcet}}
\end{center}

\title[]{}

\maketitle

\null

\vskip-60pt

\null

\begin{center}
{\large {\bf Abstract}}
\end{center}

\vskip-30pt

\null

\begin{abstract}
We investigate $L_p$-estimates for balanced averages of Fourier truncations in group algebras, in terms of differential operators acting on them. Our results extend a fundamental inequality of Naor for the hypercube (with profound consequences in metric geometry) to discrete groups. Different inequalities are established in terms of directional derivatives which are constructed via affine representations determined by the Fourier truncations. Our proofs rely on the Banach $\mathrm{X}_p$ nature of noncommutative $L_p$-spaces and dimension-free estimates for noncommutative Riesz transforms. In the particular case of free groups we use an alternative approach based on free Hilbert transforms.    
\end{abstract}

\addtolength{\parskip}{+1ex}

\vskip15pt

\section*{\bf \large Introduction}

This paper is motivated by a recent inequality due to Assaf Naor, which we now introduce. Let $\Omega_n$ be the $n$-hypercube $\Z_2 \times \Z_2 \times \cdots \times \Z_2$ equipped with its normalized counting measure. If $[n] := \{1,2,\ldots,n\}$ and we use multiplicative terminology with $\Z_2 = \{\pm 1\}$, every function $f: \Omega_n \to \C$ admits a Fourier expansion in terms of Walsh characters $W_\mathsf{A}$, which are defined by $\mathsf{A}$-products of Rademacher coordinates for any $\mathsf{A} \subset [n]$. Given a mean-zero $f$, Naor proved in \cite{N16} the following inequality for each  $p \ge 2$ and $k \in [n]$
\begin{equation} \label{Eq-Naor_p} \tag{$\mathrm{N}_p$}
\frac{1}{{{n}\choose{k}}} \sum_{\begin{subarray}{c} \mathsf{S} \hskip1pt \subset \hskip1pt [n] \\ |\mathsf{S}|=k \end{subarray}} \Big\| \sum_{\mathsf{A} \subset \mathsf{S}} \widehat{f}(\mathsf{A}) W_\mathsf{A} \Big\|_{L_p(\Omega_n)}^p \, \lesssim_p \, \frac{k}{n} \sum_{j=1}^n \|\partial_j f\|_{L_p(\Omega_n)}^p + \Big( \frac{k}{n} \Big)^{\frac{p}{2}} \|f\|_{L_p(\Omega_n)}^p.
\end{equation}
The above $\mathsf{S}$-truncations of the Walsh expansion of $f$ are conditional expectations denoted by $\mathsf{E}_{[n] \setminus \mathsf{S}}f$, while $\partial_j f$ stands for the $j$-th directional (discrete) derivative of $f$, given by $\varepsilon \mapsto f(\varepsilon) - f(\varepsilon - 2\varepsilon_j e_j)$. This inequality has groundbreaking applications in metric geometry. More precisely, it implies the quantitatively optimal form of the so-called $\mathrm{X}_p$ inequality, introduced by Naor and Schechtman in \cite{NS16}. In turn, this gives a purely metric criterion to estimate the $L_p$-distortion of a metric space $\mathrm{X}$ from below. Its metric nature is extremely useful in solving nonlinear problems around the nonembedability of $L_q$ into $L_p$ for $2 < q < p$. This includes, beyond the scope of linear $L_p$-embedding theory, the optimal $L_p$-distortion of (nonlinear) grids in $\ell_q^n$ or the critical $L_p$ snowflake exponent of $L_q$. In conclusion, Naor's differential inequality \eqref{Eq-Naor_p} and subsequent $\mathrm{X}_p$ inequalities with sharp scaling parameter are a key contribution to the Ribe program, an effort to identify which properties from the local theory of Banach spaces ultimately rely on purely metric considerations and not on the whole strength of the linear structure of the space. 

Naor's inequality \eqref{Eq-Naor_p} for functions with a linear Walsh expansion becomes a form of Rosenthal inequality for symmetrically exchangeable random variables \cite{JMST79,R70}. More precisely, let $\Pi_k$ be the space of sets $\mathsf{S} \subset [n]$ with $|\mathsf{S}|=k$ equipped with its normalized counting measure and define $\Sigma_{n,k} = \Omega_n \otimes \Pi_k$. Then, if $\widehat{f}(\mathsf{A})=0$ when $|\mathsf{A}| \neq 1$, the left hand side of \eqref{Eq-Naor_p} becomes $$\Big\| \sum_{j=1}^n \widehat{f}(\{j\}) \sigma_j \Big\|_{L_p(\Sigma_{n \hskip-0.5pt ,k})}^p \quad \mbox{with} \quad \sigma_j(\varepsilon, \mathsf{S}) = \varepsilon_j \otimes \delta_{j \in \mathsf{S}},$$ and the linear model for Naor's inequality follows from \cite{JMST79}
$$\Big\| \sum_{j=1}^n \widehat{f}(\{j\}) \sigma_j \Big\|_{L_p(\Sigma_{n \hskip-0.5pt ,k})} \asymp_p \Big( \frac{k}{n} \sum_{j=1}^n |\widehat{f}(\{j\})|^p \Big)^{\frac1p} + \Big( \frac{k}{n} \sum_{j=1}^n |\widehat{f}(\{j\})|^2 \Big)^{\frac12}.$$
Its general form \eqref{Eq-Naor_p} can be regarded as an extension for Rademacher chaos. Our primary goal in this paper is to produce similar inequalities when we replace the hypercube by other (nonnecessarily abelian) discrete groups. Fourier series with frequencies on a given discrete group $\G$ must be written in terms of its left regular representation $\lambda: \G \to \mathcal{B}(\ell_2(\G))$. The unitaries $\lambda(g)$ replace Walsh characters and we work with operators of the form $$f = \sum_{g \in \G} \widehat{f}(g) \lambda(g).$$ 
The \lq\lq quantum" probability space where we place them is the group von Neumann algebra $\V$,  precise definitions are given in the body of the paper. Understanding how to replace Rademacher chaos by some sort of \lq\lq trigonometric chaos" has to do with identifying elementary generating families. Our construction is somehow delicate and we start with a model case which originally motivated us.   

Let $\F_n = \Z * \Z * \cdots * \Z$ be the free group with $n$ generators $g_1, g_2, \ldots, g_n$. The unitaries $\lambda(g_j)$ are an archetype of Voiculescu's free random variables, which play the role of Rademacher variables above. The tensor products $\zeta_j(\mathsf{S}) = \lambda(g_j) \otimes \delta_{j \in \mathsf{S}}$ in $\Sigma_{n,k}' = \mathcal{L}(\F_n) \otimes \Pi_k$ satisfy the inequality 
$$\Big\| \sum_{j=1}^n \widehat{f}(g_j) \zeta_j \Big\|_{L_p(\Sigma_{n \hskip-0.5pt ,k}')} \asymp_p \, \Big( \frac{k}{n} \sum_{j=1}^n |\widehat{f}(g_j)|^p \Big)^{\frac1p} + \Big( \frac{k}{n} \sum_{j=1}^n |\widehat{f}(g_j)|^2 \Big)^{\frac12}.$$
The desired free form of Naor's inequality looks as follows   
\begin{equation} \label{Eq-FreeNaor_p} \tag{$\mathrm{FN}_p$}
\frac{1}{{{n}\choose{k}}} \sum_{\begin{subarray}{c} \mathsf{S} \hskip1pt \subset \hskip1pt [n] \\ |\mathsf{S}|=k \end{subarray}} \Big\| \sum_{w \in \F_\mathsf{S}} \widehat{f}(w) \lambda(w) \Big\|_p^p \, \lesssim_p \, \frac{k}{n} \sum_{j=1}^n \big\| \partial_jf \big\|_p^p + \Big( \frac{k}{n} \Big)^{\frac{p}{2}} \|f\|_p^p.
\end{equation} 
Here $\F_\mathsf{S}$ denotes the free subgroup with generators in $\mathsf{S}$ and 
$$\partial_j f= 2\pi i \Big( \sum_{w \ge g_j} \widehat{f}(w) \lambda(w) + \hskip-5pt \sum_{w \ge g_j^{-1}} \widehat{f}(w) \lambda(w) \Big),$$ where $w \ge g_j$ is used to pick those words starting with the letter $g_j$ when written in reduced form. Let us briefly comment on the two inequalities above. The first one follows from the noncommutative Burkholder/Rosenthal inequality \cite{JX03,JX08}. On the other hand, the second inequality reduces to the first one when $f$ lives in the linear span of $\lambda(g_j)$'s as a consequence of the free Khintchine inequality \cite{HP93}. It is therefore an extension of the linear model for free chaos. A  look at Naor's original inequality shows that both group elements and collections of generators (respectively denoted by $\mathsf{A}$ and $\mathsf{S}$ there) become subsets of $[n]$. This curious coincidence in the hypercube must be decoupled for other discrete groups and our inner sum in the left hand side is taken over those words $w$ with letters living in free coordinates located in $\mathsf{S}$. On the other hand, our choice for $\partial_j f$ comes from \cite{JMP18} and will be properly justified in due time. It is worth mentioning that some nonlinear extensions of the free Rosenthal inequality where investigated in \cite{JPX07} for free chaos, but none of them include a free form of Naor's inequality along the lines suggested above. 

The above reasoning settles a free model for Naor's inequality and illustrates how trigonometric chaos fits in for free groups. 
What happens if we take products of more general discrete groups? What about discrete groups lacking a product structure? Answering these questions amounts to considering Fourier truncations and somehow related differential operators over discrete groups. Other than lattices of Lie groups, discrete groups fail to admit canonical differential structures. This difficulty was successfully solved in \cite{JMP14,JMP18} with affine representations. More precisely, an orthogonal cocycle of $\G$ is a pair $(\alpha,\beta)$ given by an orthogonal action $\alpha: \G \curvearrowright \Hcal$ into some $\R$-Hilbert space together with a map $\beta: \G \to \Hcal$ satisfying the cocycle law $$\alpha_g(\beta(h)) = \beta(gh) - \beta(g).$$ The latter ensures that $g \mapsto \alpha_g(\cdot) + \beta(g)$ is an affine representation of $\G$, so that the cocycle map $\beta$ establishes a good Hilbert space lift of $\G$ and one can expect to import the differential structure of $\Hcal$. Naively, we \lq\lq identify" the unitary $\lambda(g)$ with the Euclidean character $\exp (2\pi i \langle \beta(g), \cdot \rangle)$ and define $\Hcal$-directional derivatives on $\V$ as follows for any $u \in \Hcal$
$$\partial_u(\lambda(g)) = 2\pi i \langle \beta(g), u \rangle \lambda(g) \quad \mbox{and} \quad \Delta(\lambda(g)) = -4\pi^2 \|\beta(g)\|^2 \lambda(g).$$ This strategy has been extremely useful to establish $L_p$-boundedness criteria for Fourier multipliers on group von Neumann algebras. We now introduce the right setup for the problem. Given a discrete group $\G$ equipped with an orthogonal cocycle $(\alpha,\beta)$ and a positive integer $n$, we say that $$\A = \Big\{ \B_\mathsf{S} \subset \G: \, \mathsf{S} \subset [n] \Big\}$$ is an \emph{admissible family of Fourier truncations} when we have:
\begin{itemize}
\item $\displaystyle \Big\| \sum_{g \in \B_\mathsf{S}} \widehat{f}(g) \lambda(g) \Big\|_p \le_{\mathrm{cb}} C_p \Big\| \sum_{g \in \G} \widehat{f}(g) \lambda(g) \Big\|_p$ for $p \ge 2$. 

\vskip3pt

\item Pairwise $\beta$-orthogonality: $$\Hcal = \displaystyle \bigoplus_{j=1}^n \Hcal_j \quad \mbox{with} \quad \beta(\B_\mathsf{S}), \beta(\B_\mathsf{S}^{-1})  \subset  \bigoplus_{j \in \mathsf{S}} \Hcal_j = \Hcal_\mathsf{S}.$$
\end{itemize} 
Given an orthonormal basis $(u_{j\ell})_\ell$ of $\Hcal_j$, define the \emph{$j$-th gradients} 
$$\mathrm{D}_jf = \sum_{\ell \ge 1} \partial_{u_{j\ell}}f \otimes e_{1\ell} \quad \mbox{so that} \quad |\mathrm{D}_jf| = \Big( \sum_{\ell \ge 1} |\partial_{u_{j\ell}}f|^2 \Big)^\frac12.$$

\begin{ltheorem} \label{thm:theoremA}
Let $\G$ be a discrete group equipped with an orthogonal cocycle $(\alpha,\beta)$ whose associated laplacian $\Delta$ has a positive spectral gap $\sigma > 0$. Let us consider an admissible family of Fourier truncations $\A = \{\B_\mathsf{S}: \mathsf{S} \subset [n]\}$. Then, given $p \ge 2$ and $k \in [n]$, the following inequality holds for any mean-zero $f$ $$\frac{1}{{{n}\choose{k}}} \sum_{|\mathsf{S}|=k} \Big\| \sum_{g \in \B_\mathsf{S}} \widehat{f}(g) \lambda(g) \Big\|_p^p \, \lesssim_{p,\sigma} \, \frac{k}{n} \sum_{j=1}^n \Big[  \big\|  \mathrm{D}_j(f) \big\|_p^p + \big\| \mathrm{D}_j(f^*) \big\|_p^p  \Big] + \Big( \frac{k}{n} \Big)^{\frac{p}{2}} \|f\|_p^p.$$  
\end{ltheorem}

Naor's inequality follows as a particular case of Theorem \ref{thm:theoremA} by taking $\G = \widehat{\Omega}_n$ equipped with the trivial cocycle into the $n$-dimensional space $\Hcal = \R \times \R \times \cdots \times \R$ and the truncations $\mathrm{B}_\mathsf{S} = \{\mathsf{A} \subset \mathsf{S}\} = \beta^{-1}(\Hcal_\mathsf{S})$. Recall that $\mathrm{D}_j = \partial_j$ in this case since $\dim \Hcal_j = 1$. Moreover, $|\partial_j(f)| = |\partial_j(f^*)|$ in the abelian framework of the hypercube. Two generalizations of Naor's inequality for large classes of discrete groups easily follow from Theorem \ref{thm:theoremA}: 
\begin{itemize}
\item[i)] \textbf{Direct products.} If $\G = \G_1 \times \G_2 \times \cdots \times \G_n$ is a direct product of discrete groups equipped with orthogonal cocycles $(\alpha_j,\beta_j)$, consider the product cocycle $(\alpha,\beta)$ and set $\B_\mathsf{S}$ be the subgroup of $\G$ generated by group elements whose nontrivial entries lie in $\mathsf{S}$. Then, the Fourier truncations become (completely contractive) conditional expectations and we get an admissible family of Fourier truncations. The gradients $\mathrm{D}_j$ correspond to the different factors and cocycles in the direct product above.  

\vskip3pt

\item[ii)] \textbf{Equivariant decompositions.} If $\G$ is a discrete group equipped with an orthogonal cocycle $(\alpha,\beta)$, any direct sum decomposition of the Hilbert space $\Hcal$ into $\alpha$-equivariant subspaces gives rise to an admissible family of Fourier truncations. More precisely, assume \vskip-10pt $$\hskip20pt \Hcal = \bigoplus_{j=1}^n \Hcal_j \quad \mbox{and} \quad \alpha_g(\Hcal_j) \subset \Hcal_j \mbox{ for every $(g,j) \in \G \times [n]$}.$$ Then, the family of sets \vskip-10pt $$\hskip20pt \mathrm{B}_\mathsf{S} = \beta^{-1} \Big( \bigoplus_{j \in \mathsf{S}} \Hcal_j \Big)$$ are subgroups of $\G$. In particular, the associated Fourier truncations are conditional expectations  (henceforth $L_p$-contractions) and the $\mathrm{B}_\mathsf{S}$'s satisfy pairwise $\beta$-orthogonality. This more general construction does not impose a direct product structure on the discrete group $\G$.  
\end{itemize}

Let $\A$ be an admissible family of Fourier truncations on $\G$ as defined above. Let us say that a group element $g \in \G$ is an $\A$-generator when $\beta(g) \in \Hcal_j$ for some $1 \le j \le n$.  Theorem \ref{thm:theoremA} may be regarded as a nonlinear form of an inequality for linear combinations of $\A$-generators \vskip-20pt $$f = \sum_{j=1}^n \sum_{\beta(g) \in \Hcal_j} \widehat{f}(g) \lambda(g) = \sum_{j=1}^n A_j(f).$$ This inequality controls balanced averages of $\mathsf{S}$-truncations $\sum_{j \in \mathsf{S}} A_j(f)$ in terms of $f$ and the $j$-th gradients of $A_j(f)$. This linear model seems to be new for general discrete groups/cocycles and Theorem \ref{thm:theoremA} gives a nonlinear generalization in terms of trigonometric chaos over $\A$-generators. 

Theorem \ref{thm:theoremA} does not recover the conjectured free form of Naor's inequality \eqref{Eq-FreeNaor_p}. Indeed, the free inequality relies on the standard cocycle of $\mathbb{F}_n$ associated with the word length, which yields $\Hcal \simeq \ell_2 (\mathbb{F}_n \! \setminus \! \{e\})$ and infinitely many free derivatives of the form $$\partial_u f = \sum_{w \ge u} \widehat{f}(w) \lambda(w) \quad \mbox{for any} \quad u \in \mathbb{F}_n \setminus \{e\}.$$ However, we only need to use $n$ free directional derivatives $$\partial_j = \partial_{g_j} + \partial_{g_j^{-1}} \quad  \mbox{with} \quad 1 \le j \le n$$ and these are not coupled into a family of gradients, as we do in Theorem \ref{thm:theoremA}. The key point to achieve this is the fact that free derivatives associated to free generators include all free derivatives in the sense that $$u \neq e \ \Rightarrow \ u \ge g_j \mbox{ or } u \ge g_j^{-1} \mbox{ for some $1 \le j \le n$} \ \Rightarrow \ \partial_u \circ \partial_j = \partial_j \circ \partial_u = \partial_u.$$ In general, assume that $\A = \{ \mathrm{B}_\mathsf{S}: \mathsf{S} \subset [n] \}$ is an admissible family of Fourier truncations in $\G$ with respect to $(\alpha,\beta)$. We will say that $\mathcal{J} = \{\partial_j: 1 \le j \le n\}$ is a \emph{distinguished family of derivatives} when $\partial_u \circ \partial_j = \partial_u$ for any $u \in \Hcal_j$ with $1 \le j \le n$. Throughout the paper, we shall consistently use $u$ for vectors in $\mathcal{H}$ and $j \in [n]$, so that no confusion should arise when using $\partial_u$ and $\partial_j$. The following result refines Theorem \ref{thm:theoremA} when we can find such a family.

\begin{ltheorem} \label{thm:theoremC}
Let $\G$ be a discrete group equipped with an orthogonal cocycle $(\alpha,\beta)$ and an admissible family of Fourier truncations $\A = \{\B_\mathsf{S}: \mathsf{S} \subset [n]\}$. Assume that $\mathcal{J} = \{ \partial_j: 1 \le j \le n\}$ is a distinguished family of derivatives. Then, given $p \ge 2$ and $k \in [n]$, the following inequality holds for any mean-zero $f$ $$\frac{1}{{{n}\choose{k}}} \sum_{|\mathsf{S}|=k} \Big\| \sum_{g \in \B_\mathsf{S}} \widehat{f}(g) \lambda(g) \Big\|_p^p \, \lesssim_p \, \frac{k}{n} \sum_{j=1}^n \big\| \partial_j(f) \big\|_p^p + \big\| \partial_j(f^*) \big\|_p^p + \Big( \frac{k}{n} \Big)^{\frac{p}{2}} \|f\|_p^p.$$ 
\end{ltheorem}

When the distinguished family of derivatives $\partial_j$ is a proper subset of the cocycle derivatives $\partial_u$, it turns out that Theorem \ref{thm:theoremC} gives a stronger inequality (compared to that of Theorem \ref{thm:theoremA}) at the cost of additional assumptions, which fortunately hold in several important cases considered below. Note as well that the spectral gap assumption is unnecessary under the presence of distinguished derivatives. Here are our main applications of Theorem \ref{thm:theoremC}:
\vskip-20pt \null \begin{itemize}
\item[i)] \textbf{Free chaos.} Our discussion on free derivatives illustrates how to apply Theorem \ref{thm:theoremC} to obtain an inequality which gets very close to \eqref{Eq-FreeNaor_p}. The extra term $\partial_j(f^*)$ is anyway removable due to a special property of free groups, for which word-length derivatives become free forms of directional Hilbert transforms \cite{MR17}. This \lq\lq good pathology" leads us to an even stronger inequality than the free analog of Naor's inequality \eqref{Eq-FreeNaor_p}. This could be useful in other directions of free harmonic analysis. We shall also explore the free products $\Z_{2m} * \Z_{2m} * \cdots * \Z_{2m}$.   

\vskip3pt

\item[ii)] \textbf{Continuous and discrete tori.} We also analyze $\mathbb{T}^n = \widehat{\Z}^n$ and $\Z_m^n = \widehat{\Z}_m^n$ 
equipped with different geometries. Theorem \ref{thm:theoremC} is applicable for the Cayley graph metric and the resulting inequality improves the one coming from the Euclidean metric. These forms of Naor's inequality can be regarded as refinements of the classical Poincar\'e inequality.

\vskip3pt

\item[iii)] \textbf{Infinite Coxeter groups.} Any group presented by $$\hskip32pt \G \, = \, \Big\langle g_1, g_2, \ldots, g_n \, \big| \, (g_jg_k)^{s_{jk}}= e \Big\rangle$$ with $s_{jj}=1$ and $s_{jk} \ge 2$ for $j \neq k$ is called a Coxeter group. Bo\.zejko proved in \cite{Bz} that the word length is conditionally negative for any infinite Coxeter group. The Cayley graph of these groups is more involved and we will not construct here a natural ONB for the cocycle, we invite the reader to do it and to derive inequalities in the lines of Theorems \ref{thm:theoremA} and \ref{thm:theoremC}.  
\end{itemize}

Our proof of Theorems \ref{thm:theoremA} and \ref{thm:theoremC} streamlines Naor's original argument. The key point in this general setting is to identify the right notions, such as admissible families of Fourier truncations or distinguished families of derivatives. Once this is done, the proof heavily relies on dimension-free estimates for noncommutative Riesz transforms \cite{JMP18} in the same way Naor's inequality did in terms of Lust-Piquard results \cite{LP98}. Another crucial point in our argument is the Banach $\mathrm{X}_p$ nature of noncommutative $L_p$-spaces. Generalizing previous work of Naor and Schechtman \cite[Theorem 7.1]{NS16}, we shall establish it with a much simpler argument based on Junge/Xu's noncommutative Burkholder and Rosenthal inequalities \cite{JX03,JX08}. Of course, one could expect that Theorems \ref{thm:theoremA} and \ref{thm:theoremC} may lead to noncommutative $\mathrm{X}_p$-type inequalities, very much like in \cite{N16}. We have obtained some inequalities in this direction \cite{CCP22p2}. Our hope was to deduce nontrivial bounds for $L_p$-distortions of Schatten $q$-classes or other noncommutative $L_q$-spaces. Unfortunately, our efforts so far have not been fruitful in this direction. 


\section{\bf Trigonometric chaos} \label{sec1}

\subsection{Harmonic analysis on discrete groups} \label{ssec11} Let $\G$ be a discrete group. The left regular representation of $\G$ on $\ell_2(\G)$ is the unitary representation determined by
\begin{equation*}
    [\lambda(g)\varphi](h) = \varphi(g^{-1}h), \;\; g,h \in \G, \; \varphi \in \ell_2(\G).
\end{equation*}
The group von Neumann algebra of $\G$ is denoted by $\mathcal{L}(\G)$. It is the weak operator closure of the linear span of $\{\lambda(g)\}_{g \in \G}$ in $\mathcal{B}(\ell_2(\G))$. Its canonical trace $\tau$ is linearly determined by $\tau(\lambda(g)) = \langle \lambda(g)\delta_e, \delta_e \rangle_{\ell_2(\G)} = \delta_{g=e}$. Every element $f\in \mathcal{L}(\G)$ admits a Fourier series
$$
f = \suma{g \in \G} \widehat{f}(g) \lambda(g) \quad \mbox{where} \quad \widehat{f}(g) = \tau(\lambda(g)^* f). 
$$
This shows that $\tau(f)=\widehat{f}(e)$. For $1 \leq p < \infty$, we denote by $L_p(\mathcal{L}(\G))$ the associated noncommutative $L_p$ space. An orthogonal cocycle for $\G$ is a triple $(\mathcal{H}, \alpha, \beta)$ given by a real Hilbert space $\mathcal{H}$, an orthogonal action $\alpha : \G \rightarrow \mathcal{O}(\mathcal{H})$, and a map $\beta : \G \rightarrow \mathcal{H}$ satisfying the cocycle law
$$
\alpha_g(\beta(h)) = \beta(gh) - \beta(g).
$$
Orthogonal cocycles are in one-to-one correspondence with length functions. We say that a map $\psi : \G \rightarrow \R_+$ is a length function if it vanishes at the identity $e$, it is symmetric $\psi(g)=\psi(g^{-1})$, and it is conditionally negative
$$
\sum_{g \in \G} a_g = 0 \quad \Rightarrow \quad \sum_{g,h \in \G} \overline{a}_g a_h \psi(g^{-1}h) \leq 0
$$
for any finitely supported family $\{a_g\}_{g \in \G}$. Given a cocycle $(\mathcal{H},\alpha,\beta)$, the function $\psi(g) = \|\beta(g)\|_{\mathcal{H}}^2$ is a length function. On the other hand, any length function $\psi$ determines a Gromov form $\langle \cdot, \cdot \rangle_{\psi}$, a semidefinite positive form on the group algebra $\R[\G]$ given by
$$
\langle \delta_g, \delta_h \rangle_{\psi} = \frac{\psi(g)+\psi(h)-\psi(g^{-1}h)}{2}.
$$
Then, the Hilbert completion $\mathcal{H}$ of $(\R[\G]/\mathrm{Ker(\langle \cdot, \cdot \rangle_{\psi})},\langle \cdot, \cdot \rangle_\psi)$, together with the map $\beta: g \mapsto \delta_g + \mathrm{Ker}(\langle \cdot, \cdot \rangle_{\psi})$, and the orthogonal action $\alpha_{g}(\delta_h) = \delta_{gh} - \delta_{h} + \mathrm{Ker}(\langle \cdot, \cdot \rangle_{\psi})$ form a  cocycle. Moreover, Schoenberg's theorem \cite{S38} claims that $\psi : \G \rightarrow \R_+$ is a length function iff the maps $S_{t} : \lambda(g) \mapsto e^{-t\psi(g)} \lambda(g)$ form a Markov semigroup $(S_{t})_{t \geq 0}$ on $\mathcal{L}(\G)$, see \cite{JMP14,JMP18}. In this case $(S_{t})_{t \geq 0}$ admits an infinitesimal generator
$$
-\laplacian := \lim_{t \rightarrow 0^+} \frac{S_{t} - \mathrm{id}_{\mathcal{L}(\G)}}{t} \quad \mathrm{\ so \, that \ } \quad S_{t} = \exp(-t \laplacian).
$$
As is standard, we shall call the generator $\laplacian$ the $\psi$-Laplacian on $\G$. Since we have $\laplacian(\lambda(g)) = \psi(g) \lambda(g)$ for $g\in \G$, it turns out that $\laplacian$ is an unbounded Fourier multiplier whose fractional powers can be defined by 
$$
\laplacian^{\gamma} f := \suma{g \in \G} \psi(g)^{\gamma} {f}(g) \lambda(g).
$$
Let $(\mathcal{H}, \alpha, \beta)$ be the orthogonal cocycle naturally associated to the length function $\psi: \G \to \R_+$ as explained above. Given an orthonormal basis $\{e_\ell\}_{\ell \geq 1}$ of $\mathcal{H}$, we consider the corresponding directional derivatives as follows
$$
\partial_{e_\ell} \lambda(g) := 2\pi i  \langle \beta(g), e_\ell \rangle_{\psi} \lambda(g) \quad \mbox{so that} \quad -4\pi^2 \laplacian = \sum_{\ell \ge 1} \partial_{e_\ell}^2.
$$
The corresponding Riesz transforms associated to $\psi$ are then defined as
\begin{align*}
	R_\ell f=R_{e_\ell}f := \partial_{e_\ell} \laplacian^{-\frac12} f = 2\pi i \suma{g \in \G} \frac{\langle \beta(g), e_\ell \rangle_{\psi}}{\sqrt{\psi(g)}} \, \widehat{f}(g) \lambda(g).
\end{align*} 
Riesz transforms act on elements of $L_p(\mathcal{L}(\G))$ with null Fourier coefficients on the kernel of $\beta$. More precisely, maps $R_\ell$ are well-defined over the mean-zero subspaces
$$
	L_p^\circ(\mathcal{L}(\G)) = \Big\{ f \in L_p(\mathcal{L}(\G)) : \widehat{f}(g)= 0 \mathrm{\ if \ } \beta(g)=0 \Big\}.
$$
Dimension-free estimates for noncommutative Riesz transforms were studied in \cite{JMP18}.

\begin{theorem}[Theorem A1 -- \cite{JMP18}] \label{JMPDimFree}
If $2 \leq p < \infty$ and $f\in L_p^\circ(\mathcal{L}(\G))$
$$
\|f\|_p \ \asymp_p \ \max\Bigg\{ \Big\| \Big(\suma{\ell \geq 1} |R_{\ell}(f)|^2 \Big)^{\frac12} \Big\|_p, \Big\| \Big(\suma{\ell \geq 1} |R_{\ell}(f^*)|^2 \Big)^{\frac12} \Big\|_p \Bigg\}.
$$
\end{theorem}
\noindent Finally, our Fourier truncations will be written in the form
$$
\cexp_{[n]\setminus \Ss} f = \sum_{g \in \B_\Ss} \widehat{f}(g)\lambda(g) \quad \mbox{with} \quad \Ss \subset [n]. 
$$
When $\B_\Ss$ is a subgroup of $\G$, $\cexp_{[n]\setminus \Ss}$ is a ($L_p$-contractive) conditional expectation. 

\subsection{Noncommutative $L_p$-spaces are Banach $\mathrm{X}_p$ spaces}

Linear forms of $\mathrm{X}_p$ inequalities are vector-valued extensions of Rosenthal inequality for symmetrically exchangeable random variables \cite{JMST79}. In \cite[Theorem 7.1]{NS16} Naor and Schechtman proved such inequalities for Schatten $p$-classes. A noncommutative Burkholder martingale inequality for the conditioned square function \cite{JX03} led Junge and Xu to obtain noncommutative Rosenthal inequalities for symmetric variables in \cite{JX08}. We use this result below to establish the Banach $\mathrm{X}_p$ nature of arbitrary noncommutative $L_p$-spaces. Naor/Schechtman's argument can be extended to work as well for other noncommutative $L_p$-spaces, but our argument below is much shorter. 

\begin{theorem} \label{thm:linearXp}
Let $(\mathcal{M},\tau)$ be a von Neumann algebra equipped with a normal semifinite faithful trace. Then, if $\mathbb{E}$ denotes the expectation over independently equidistributed random signs $\varepsilon = (\varepsilon_1, \varepsilon_2, \ldots, \varepsilon_n)$ and $x_j \in L_p(\mathcal{M})$, the following inequality holds for any $p \ge 2$ and $k \in [n]$
$$
\frac{1}{{n \choose k}} \suma{\substack{\Ss \subseteq [n] \\ |\Ss|=k}} \mathbb{E} \Big\| \suma{j \in \Ss} \varepsilon_j x_j \Big\|^p_{L_p(\mathcal{M})} \lesssim_p \frac{k}{n} \Suma{j=1}{n} \|x_j\|^p_{L_p(\mathcal{M})} + \Big(\frac{k}{n}\Big)^{\frac{p}{2}} \mathbb{E} \Big\| \Suma{j=1}{n} \varepsilon_j x_j \Big\|_{L_p(\mathcal{M})}^p.
$$
\end{theorem}

\dem Define random variables $\sigma_j \in \Sigma_{n,k} = \Omega_n \otimes \Pi_k$ as defined in the Introduction by $\sigma_j(\varepsilon,\Ss) = \varepsilon_{j} \otimes \delta_{j \in \Ss}$ for $1 \le j \le n$ and $\Ss\subset [n]$. We claim that, for any choice of signs $\delta_j = \pm 1$ and any permutation $\pi$ of $[n]$, it holds
\begin{equation*} 
\mathrm{A} := \Big\| \Suma{j=1}{n} \delta_j \sigma_{\pi(j)} \otimes x_j \Big\|_{L_p(\Sigma_{n,k} \overline{\otimes} \mathcal{M})} \lesssim \Big\|\Suma{j=1}{n} \sigma_j \otimes x_j \Big\|_{L_p(\Sigma_{n,k} \overline{\otimes} \mathcal{M})} =: \mathrm{B}.
\end{equation*} 
Indeed, applying noncommutative Khintchine's inequality \cite{LP86} twice
\begin{eqnarray*}
\mathrm{A}^p \!\!\! & = & \!\!\! \frac{1}{{n \choose k}} \suma{\substack{\Ss\subseteq [n] \\ |\Ss|=k}} \Big\| \suma{\pi(j) \in \Ss} \varepsilon_{\pi(j)} \otimes \delta_j x_j \Big\|^p_{L_p(\Omega_n;L_p(\mathcal{M}))} \\
\!\!\! & \asymp_p & \!\!\! \frac{1}{{n \choose k}} \suma{\substack{\Ss \subseteq [n] \\ |\Ss|=k}} \max \Bigg\{ \Big\| \Big( \suma{\pi(j) \in \Ss} x_j^* x_j \big)^\frac12 \Big\|_{L_p(\mathcal{M})}, \Big\| \Big( \suma{\pi(j) \in \Ss} x_j x_j^* \Big)^\frac12 \Big\|_{L_p(\mathcal{M})} \Bigg\} \\
\!\!\! & = & \!\!\! \frac{1}{{n \choose k}} \suma{\substack{\Ss \subseteq [n] \\ |\Ss|=k}} \max \Bigg\{ \Big\| \Big( \suma{j \in \Ss} x_j^* x_j \big)^\frac12 \Big\|_{L_p(\mathcal{M})}, \Big\| \Big( \suma{j \in \Ss} x_j x_j^* \big)^\frac12 \Big\|_{L_p(\mathcal{M})} \Bigg\} \ \asymp_p \ \mathrm{B}^p.
\end{eqnarray*}
Hence, we can apply \cite[Corollary 6.6]{JX08} to get
$$
\mathrm{B}^p \lesssim_p \frac{1}{n} \Suma{j,j'=1}{n} \|\sigma_j\|_p^p \|x_{j'}\|_p^p
 + \Big(\frac{1}{n}\Big)^{\frac{p}{2}} \Big\| \Big(\Suma{j=1}{n} x_j^* x_j + x_j x_j^*\Big)^\frac12 \Big\|_p^p \Big\|\Big(\Suma{j=1}{n} \sigma_j^2\Big)^\frac12 \Big\|_p^p.
$$
Now, we have
\begin{eqnarray*}
\|\sigma_j\|_{L_p(\Sigma_{n,k})}^p \!\!\! & = & \!\!\! \frac{1}{{n \choose k}} \suma{\substack{\Ss \subseteq [n] \\ |\Ss|=k}} \delta_{j \in \Ss} = \frac{k}{n}, \\ \Big\|\Big( \Suma{j=1}{n} \sigma_j^2\Big)^\frac12 \Big\|_{L_p(\Sigma_{n,k})}^p \!\!\! & = & \!\!\! \frac{1}{{n \choose k}} \suma{\substack{\Ss \subseteq [n] \\ |\Ss|=k}} \Big( \Suma{j=1}{n} \delta_{j \in \Ss} \Big)^{\frac{p}{2}} =  k^{\frac{p}{2}}.
\end{eqnarray*}
Therefore, we get
\begin{eqnarray*}
\mathrm{B} \!\!\! & \lesssim_p & \!\!\! \frac{k}{n} \Suma{j=1}{n} \|x_j\|_{L_p(\mathcal{M})}^p + \Big(\frac{k}{n}\Big)^{\frac{p}{2}} \Big\|\Big( \Suma{j=1}{n} x_j^* x_j + x_j x_j^* \Big)^\frac12 \Big\|_{L_p(\mathcal{M})}^p \\ \!\!\! & \asymp_p & \!\!\! \frac{k}{n} \Suma{j=1}{n} \|x_j\|_{L_p(\mathcal{M})}^p + \Big(\frac{k}{n}\Big)^{\frac{p}{2}} \mathbb{E} \Big\| \Suma{j=1}{n} \varepsilon_j x_j \Big\|_{L_p(\mathcal{M})}^p, 
\end{eqnarray*} 
applying once again noncommutative Khintchine's inequality. This proves the result since the random variables $\sigma_j$ are chosen so that $\mathrm{B}$ equals the left hand side in the inequality of the statement. \fin

\begin{remark} 
\emph{Theorem \ref{thm:linearXp} says that $L_p(\mathcal{M})$ is an Banach $\mathrm{X}_p$ space. The conclusion also holds in the completely bounded setting since the constants that appear in the inequality of the statement do not depend on the von Neumann algebra $\mathcal{M}$.}
\end{remark}

\subsection{Proof of Theorem \ref{thm:theoremA}}

According to Theorem \ref{JMPDimFree}
\begin{eqnarray*}
\frac{1}{{{n}\choose{k}}} \sum_{\substack{\Ss \subseteq [n] \\ |\Ss|=k}} \Big\| \sum_{g \in \B_\mathsf{S}} \widehat{f}(g) \lambda(g) \Big\|_p^p \!\!\!\! & = & \!\!\!\! \frac{1}{{n \choose k}} \suma{\substack{\Ss \subseteq [n] \\ |\Ss|=k}} \big\| \cexp_{[n]\setminus \Ss}f \big\|_p^p \\ \!\!\!\! & \asymp_p & \!\!\!\!  \frac{1}{{n \choose k}} \suma{\substack{\Ss \subseteq [n] \\ |\Ss|=k}} \Big\| \Big( \sum_{\substack{j \in [n] \\ \ell \ge 1}} |R_{j\ell}(\cexp_{[n]\setminus \Ss}f)|^2\Big)^{\frac12} \Big\|_p^p \nonumber \\ \!\!\!\! & + & \!\!\!\! \frac{1}{{n \choose k}} \suma{\substack{\Ss \subseteq [n] \\ |\Ss|=k}} \Big\| \Big(\sum_{\substack{j \in [n] \\ \ell \ge 1}} |R_{j\ell}((\cexp_{[n]\setminus \Ss}f)^*)|^2\Big)^{\frac12} \Big\|_p^p =: \mathrm{A+B}, 
\end{eqnarray*}
where $R_{j\ell} := R_{u_{j\ell}}$ and $\{u_{j\ell}: j \in [n], \, \ell \ge 1\}$ is the orthonormal basis of $\mathcal{H}$ considered before the statement of Theorem \ref{thm:theoremA}. Since $\beta(\B_\Ss) \subset \mathcal{H}_\Ss$, we observe that $\langle \beta(g), u_{j\ell} \rangle_\psi = 0$ whenever $g \in \B_\Ss$ and $j \notin \Ss$. Moreover, Fourier truncations commute with Riesz transforms |as both are Fourier multipliers| and we deduce $$R_{jk} \circ \cexp_{[n]\setminus \Ss} = \delta_{j \in \Ss} \ \cexp_{[n]\setminus \Ss} \circ R_{jk}.$$ Using the complete $L_p$-boundedness of our Fourier truncations, we get 
\begin{eqnarray*}
\mathrm{A} \!\!\! & \lesssim_p & \!\!\! \frac{1}{{n \choose k}} \suma{\substack{\Ss \subseteq [n] \\ |\Ss|=k}} \Big\| \Big( \sum_{j\in \Ss} \sum_{\ell\geq 1} |R_{u_{j\ell}} f|^2 \Big)^{\frac12} \Big\|_p^p \\ \!\!\! & \lesssim_p & \!\!\! \frac{1}{{n \choose k}} \suma{\substack{\Ss \subseteq [n] \\ |\Ss|=k}} \mathbb{E} \Big\| \sum_{j\in \Ss} \varepsilon_j \Big[ \sum_{\ell\geq 1} |R_{u_{j\ell}} f|^2 \Big]^{\frac12} \Big\|_p^p =: \mathrm{A}'.
\end{eqnarray*}
The last inequality follows from the noncommutative Khintchine inequality \cite{LP86} applied to independent equidistributed signs $\varepsilon_j = \pm 1$. Next, we use the Banach $\mathrm{X_p}$ nature of noncommutative $L_p$-spaces. More precisely, applying Theorem \ref{thm:linearXp} we get
$$\mathrm{A}' \lesssim_p \frac{k}{n} \Suma{j = 1}{n} \Big\| \Big( \sum_{\ell \geq 1} |R_{u_{j\ell}} f|^2 \Big)^{\frac12} \Big\|_p^p + \Big(\frac{k}{n} \Big)^{\frac{p}{2}} \, \mathbb{E} \Big\| \Suma{j=1}{n} \varepsilon_j \Big( \sum_{\ell\geq 1} |R_{u_{j\ell}} f|^2  \Big)^{\frac12} \Big\|_p^p = \mathrm{A}_1' + \mathrm{A}_2'. $$
Since $R_{u_{j\ell}} = \partial_{u_{j\ell}} \Delta^{-\frac12} = \Delta^{-\frac12} \partial_{u_{j\ell}}$, \cite[Proposition 1.1.5]{JM10} yields
\begin{eqnarray*}
\mathrm{A}_1' \!\!\! & = & \!\!\! \frac{k}{n} \Suma{j=1}{n} \Big\| \sum_{\ell \geq 1} R_{u_{j\ell}} f \otimes e_{\ell,1} \Big\|_{S_p[L_p(\mathcal{L}(\G))]}^p \\ \!\!\! & \lesssim_{p,\sigma} & \!\!\! \frac{k}{n} \Suma{j=1}{n} \Big\| \sum_{\ell \geq 1} \, \partial_{u_{j\ell}} f \otimes e_{\ell,1} \Big\|_{S_p[L_p(\mathcal{L}(\G))]}^p = \frac{k}{n} \Suma{j = 1}{n} \big\| \gradient_j (f) \big\|_p^p.
\end{eqnarray*}		
Moreover, noncommutative Khintchine inequality and Theorem \ref{JMPDimFree} give
$$
\mathrm{A}_2' \, \lesssim_p \, \Big( \frac{k}{n}\Big)^{\frac{p}{2}} \Big\| \Big( \Suma{j=1}{n} \sum_{\ell \geq 1} |R_{u_{j\ell}} f|^2 \Big)^\frac12 \Big\|_p^p \, \lesssim_p \, \Big( \frac{k}{n}\Big)^{\frac{p}{2}} \|f\|_p^p.
$$
Therefore, the term $\mathrm{A}$ satisfies the expected estimate and it remains to justify the assertion for $\mathrm{B}$. We now analyze the behavior of our Fourier truncations under adjoints. Observe that
$$
(\cexp_{[n]\setminus \Ss}f)^* = \sum_{g \in \B_\Ss} \overline{\widehat{f}(g)}\lambda(g^{-1}) =: \cexp_{[n]\setminus \Ss}'(f^*).
$$
In particular, since $\cexp_{[n]\setminus \Ss}'$ commutes with $R_{j\ell}$ 
$$
R_{j\ell} \big( (\cexp_{[n]\setminus \Ss}f)^* \big) = \cexp_{[n]\setminus \Ss}' \big( R_{j \ell}(f^*) \big) = \cexp_{[n]\setminus \Ss}\big( R_{j \ell}(f^*)^* \big)^*.
$$
At this point is where we need the condition $\beta(\B_\Ss^{-1}) \subset \mathcal{H}_\Ss$, to make sure that the above terms vanish when $j \notin \Ss$ since we find the inner products $\langle \beta(g^{-1}), u_{j\ell} \rangle_\psi$ for $g \in \B_\Ss$. Thus, we obtain
\begin{eqnarray*}
\Big\| \Big(\sum_{\substack{j \in [n] \\ \ell \ge 1}} |R_{j\ell}((\cexp_{[n]\setminus \Ss}f)^*)|^2\Big)^{\frac12} \Big\|_p \!\!\! & = & \!\!\! \Big\| \sum_{j \in [n]} \sum_{\ell \ge 1} \cexp_{[n]\setminus \Ss}\big( R_{j \ell}(f^*)^* \big) \otimes e_{1,(j,\ell)} \Big\|_p \\ [-5pt] \!\!\! & \lesssim_p & \!\!\! \Big\| \sum_{j \in \Ss} \sum_{\ell \ge 1} R_{j \ell}(f^*) \otimes e_{(j,\ell),1} \Big\|_p \\ \!\!\! & = & \!\!\! \Big\| \Big( \sum_{j \in \Ss} \sum_{\ell \ge 1} | R_{j \ell}(f^*) |^2 \Big)^\frac12 \Big\|_p.
\end{eqnarray*}
Therefore, we may follow the above argument for $\mathrm{A}$ just replacing $f$ by $f^*$. \fin

\begin{remark} \label{Rem-ThmB'}
\emph{A careful reading of the proof of Theorems \ref{thm:theoremA} and \ref{thm:theoremC} shows that we may use different Hilbert space decompositions $\mathcal{H} = \oplus_j \mathcal{H}_j = \oplus_j \mathcal{K}_j$ for $\mathsf{B}_\Ss$ and its inverse |with $\beta(\mathsf{B}_\Ss) \subset \mathcal{H}_\Ss$ and $\beta(\mathsf{B}_\Ss^{-1}) \subset  \mathcal{K}_\Ss$| as long as we can find an orthonormal basis $\{u_\ell: \ell \ge 1\}$ of $\mathcal{H}$ satisfying that 
\begin{equation} \label{eq:CompatibleONB}
\forall \ \ell \ge 1 \ \ \exists \ j_1,j_2\in [n] \ \ \mbox{such that} \ \ u_\ell \in \mathcal{H}_{j_1} \cap \mathcal{K}_{j_2}.
\end{equation}
More precisely, under these more flexible assumption we get}
$$
\frac{1}{{{n}\choose{k}}} \sum_{|\mathsf{S}|=k} \Big\| \sum_{g \in \B_\mathsf{S}} \widehat{f}(g) \lambda(g) \Big\|_p^p \, \lesssim_{p,\sigma} \, \frac{k}{n} \sum_{j=1}^n \Big[\big\| \mathrm{D}_j(f) \big\|_p^p + \big\| \mathrm{D}_j^\dag(f^*) \big\|_p^p\Big] + \Big( \frac{k}{n} \Big)^{\frac{p}{2}} \|f\|_p^p,
$$
\emph{where $\gradient_j^\dag := \sum_{u_\ell \in \mathcal{K}_j} \partial_{u_{\ell}}(\cdot) \otimes e_{1\ell}$ is the gradient over the basis vectors living in $\mathcal{K}_j$.}
\end{remark}

\begin{remark}
\emph{The constant depending on $\sigma$ in Theorem \ref{thm:theoremA} grows as $\sigma^{-\frac{p}{2}}$. One can also track the dependence on $p$ of the constant. Using free generators in place of random signs |Theorem \ref{thm:linearXp} holds as well| we keep constants uniformly bounded replacing noncommutative by free Khintchine inequalities \cite{HP93}. The constants in Theorem \ref{JMPDimFree} are bounded by $p^{3/2}$, but it is still open whether this is optimal.}
\end{remark}

\subsection{Proof of Theorem \ref{thm:theoremC}} Again Theorem \ref{JMPDimFree} gives 
$$
\frac{1}{{n \choose k}} \suma{\substack{\Ss \subseteq [n] \\ |\Ss|=k}} \|\cexp_{[n]\setminus \Ss}f\|_p^p \, \asymp_p \, \mathrm{A+B}
$$
as in the proof of Theorem \ref{thm:theoremA}. Following our argument there, we use our estimate $\mathrm{A} \lesssim_p \mathrm{A}_1' + \mathrm{A}_2'$ and we bound $\mathrm{A}_2'$ in the same way. To estimate $\mathrm{A}_1'$ we use our distinguished family of derivatives and Theorem \ref{JMPDimFree} to deduce
\begin{align*}
	\mathrm{A}_1' = \frac{k}{n} \Suma{j = 1}{n} \Big\| \Big( \sum_{\ell \geq 1} |R_{u_{j\ell}} \partial_j f|^2 \Big)^{\frac12} \Big\|_p^p 
	&\lesssim_p \frac{k}{n} \Suma{j = 1}{n} \big\| \partial_j f \big\|_p^p. 	
\end{align*}
The estimate for $\mathrm{B}$ then follows by the same considerations as in Theorem \ref{thm:theoremA}. \fin


\section{\bf Applications to abelian groups} \label{sec2}

We now focus our attention on concrete realizations of Theorems \ref{thm:theoremA} and \ref{thm:theoremC} for certain commutative group algebras. In all the cases in this section, we choose $\cexp_{[n]\setminus \Ss}$ of the form
$$
\cexp_{[n]\setminus \Ss} f = \suma{g \in \B_\Ss} \widehat{f}(g) \lambda(g) \quad \mbox{for some subgroup $\B_\Ss$ of $\G$}. 
$$
Due to that fact, we know that they are conditional expectations, and therefore completely contractive maps. This allows us to safely apply Theorems \ref{thm:theoremA} and \ref{thm:theoremC} without checking that hypothesis. We will give the details for the cases $\G=\mathbb{Z}^n$ and $\G=\mathbb{Z}_{2m}^n$, yielding inequalities in $L_p(\mathbb{T}^n)$ and $L_p(\mathbb{Z}_{2m}^n)$, respectively.

\subsection{Classical tori} Define 
\begin{eqnarray*} 
\psi_1(g) \!\!\! & = & \!\!\! |g_1| + \ldots + |g_n|, \\
\psi_2(g) \!\!\! & = & \!\!\! g_1^2 + g_2^2 + \ldots + g_n^2,
\end{eqnarray*}
with $g =(g_1,g_2,\ldots,g_n) \in \Z^n$. Both functions are symmetric and vanish at $0$. Moreover, conditional negativity follows easily. In the case of $\psi_1$, it suffices to check it for each summand $|g_j|$ which are conditionally negative from subordination with respect to $g_j^2$. These functions are denoted as the word and the Euclidean length respectively. We analyze balanced Fourier truncations using both geometries.  

\noindent A) The Euclidean length. The length $\psi_2$ induces the standard cocycle $(\mathcal{H},\alpha,\beta)$ where $\mathcal{H} = \R^n$ with the usual Euclidean inner product, the trivial action and the canonical inclusion $\beta = \mathrm{id}$. We use the standard decomposition $\mathcal{H} = \bigoplus_j \mathcal{H}_j$ with $\mathcal{H}_j = \R e_j$ the subspace generated by the $j$-th element of the canonical basis. Therefore, given $\Ss\subset [n]$, denote by $\mathbb{Z}^\Ss$ the subgroup of elements with vanishing entries outside $\Ss$ and consider the truncations 
\begin{equation*}
\cexp_{[n]\setminus \Ss} f (x) \, = \, \suma{g \in \Z^\Ss}^{\null} \widehat{f}(g) e^{2\pi i \langle x, g \rangle} \quad \mbox{for any} \quad f \in L_p(\mathbb{T}^n) \simeq L_p(\mathcal{L}(\Z^n)),
\end{equation*}
where $e^{2 \pi i \langle \cdot , g \rangle} \mapsto \lambda(g)$ defines a trace-preserving $*$-homomorphism. The cocycle derivatives coincide in this case with the classical ones $\partial_{e_j} \lambda(g) = 2\pi i g_j \lambda(g)$ and the infinitesimal generator $\laplacian$ is the usual Laplacian (up to a universal constant) with spectral gap $1$. Then, Theorem \ref{thm:theoremA} yields 
\begin{equation}\label{eq:ZntheoremA}
\frac{1}{{n \choose k}} \suma{\substack{\Ss \subseteq [n] \\ |\Ss|=k}} \Big\| \suma{g \in \Z^\Ss}^{\null} \widehat{f}(g) e^{2\pi i \langle \cdot, g \rangle} \Big\|^p_{L_p(\mathbb{T}^n)} \lesssim_p \frac{k}{n} \Suma{j = 1}{n} \| \partial_{e_j} f \|_{L_p(\mathbb{T}^n)}^p + \Big(\frac{k}{n}\Big)^{\frac{p}{2}} \|f\|^p_{L_p(\mathbb{T}^n)}
\end{equation}
for any mean-zero $f\in L_p(\mathbb{T}^n)$. This seems to be the most natural generalization of Naor's inequality for classical tori, but it is not the most efficient. Indeed, using the same Hilbert space decomposition as above, one can consider the alternative absorbent derivatives $\partial_j \lambda(g) = \delta_{g_j \neq 0} \lambda(g)$, which satisfy $\partial_{e_j} \circ \partial_j = \partial_{e_j}$. In particular Theorem \ref{thm:theoremC} yields
\begin{equation}\label{eq:ZntheoremC}
\frac{1}{{n \choose k}} \suma{\substack{\Ss \subseteq [n] \\ |\Ss|=k}} \Big\| \suma{g \in \Z^\Ss}^{\null} \widehat{f}(g) e^{2\pi i \langle \cdot, g \rangle} \Big\|^p_{L_p(\mathbb{T}^n)} \lesssim_p \frac{k}{n} \Suma{j = 1}{n} \| \partial_j f \|_{L_p(\mathbb{T}^n)}^p + \Big(\frac{k}{n}\Big)^{\frac{p}{2}} \|f\|^p_{L_p(\mathbb{T}^n)}.
\end{equation}
This is a stronger inequality since $$\| \partial_j f \|_{L_p(\mathbb{T}^n)} = \frac{1}{2\pi} \Big\| \sum_{g_j \neq 0} \frac{1}{g_j} \widehat{\partial_{e_j} f}(g) e^{2\pi i \langle \cdot, g \rangle} \Big\|_{L_p(\mathbb{T}^n)} \le C_p \| \partial_{e_j} f \|_{L_p(\mathbb{T}^n)}.$$ Indeed, the symbol $m(g) = 1/g_j$ defines an $L_p$-bounded multiplier as a consequence of K. de Leeuw restriction theorem and H\"ormander-Mikhlin multiplier theorem \cite{dL,Ho,Mi}. As we shall see \eqref{eq:ZntheoremC} naturally appears using the word length. 

\begin{remark}
\emph{Consider $f: \mathbb{T}^n \to \C$ with $$f(x) = \sum_{g \in \Z^n} \widehat{f}(g) e^{2\pi i \langle x,g \rangle} \quad \mbox{and} \quad \widehat{f}(0)=0.$$ Given $\mathsf{S} \subset [n]$, the classical Poincar\'e inequality gives}
\begin{eqnarray*}
\frac{1}{{{n}\choose{k}}} \sum_{\begin{subarray}{c} \mathsf{S} \hskip1pt \subset \hskip1pt [n] \\ |\mathsf{S}|=k \end{subarray}} \Big\| \underbrace{\sum_{g \in \Z^{\mathsf{S}}} \widehat{f}(k) e^{2\pi i \langle \cdot, g \rangle}}_{f_\mathsf{S}} \Big\|_p^p \!\!\! & \le & \!\!\! \frac{1}{{{n}\choose{k}}} \sum_{\begin{subarray}{c} \mathsf{S} \hskip1pt \subset \hskip1pt [n] \\ |\mathsf{S}|=k \end{subarray}} \big\| \nabla f_\mathsf{S} \big\|_p^p \\ [-5pt] \!\!\! & \asymp & \!\!\! \frac{1}{{{n}\choose{k}}} \sum_{\begin{subarray}{c} \mathsf{S} \hskip1pt \subset \hskip1pt [n] \\ |\mathsf{S}|=k \end{subarray}} \Big\| \sum_{j \in \mathsf{S}} \varepsilon_j \partial_{e_j} f_\mathsf{S} \Big\|_p^p \\ \!\!\! & = & \!\!\! \frac{1}{{{n}\choose{k}}} \sum_{\begin{subarray}{c} \mathsf{S} \hskip1pt \subset \hskip1pt [n] \\ |\mathsf{S}|=k \end{subarray}} \Big\| \Big[ \sum_{j \in \mathsf{S}} \varepsilon_j \partial_{e_j} f \Big]_\mathsf{S} \Big\|_p^p \\ \!\!\! & \le & \!\!\! \frac{1}{{{n}\choose{k}}} \sum_{\begin{subarray}{c} \mathsf{S} \hskip1pt \subset \hskip1pt [n] \\ |\mathsf{S}|=k \end{subarray}} \Big\| \sum_{j \in \mathsf{S}} \varepsilon_j \partial_{e_j} f \Big\|_p^p = \Big\| \sum_{j=1}^n \sigma_j \partial_{e_j}f \Big\|^p_p 
\end{eqnarray*}
\emph{for $\sigma_j(\varepsilon, \mathsf{S}) = \varepsilon_j \otimes \delta_{j \in \mathsf{S}}$ as in the Introduction. Applying  \cite{JMST79} gives
$$\frac{1}{{{n}\choose{k}}} \sum_{\begin{subarray}{c} \mathsf{S} \hskip1pt \subset \hskip1pt [n] \\ |\mathsf{S}|=k \end{subarray}} \Big\| \sum_{g \in \Z^\mathsf{S}} \widehat{f}(g) e^{2\pi i \langle \cdot, g \rangle} \Big\|_p^p \, \lesssim_p \, \frac{k}{n} \sum_{j=1}^n \|\partial_{e_j} f\|_p^p + \Big( \frac{k}{n} \Big)^{\frac{p}{2}} \big\| \nabla f \big\|_p^p.$$ Inequalities \eqref{eq:ZntheoremA} and \eqref{eq:ZntheoremC} improve the above inequality replacing $\nabla f$ by $f$.}
\end{remark}

\noindent B) The word length. Let us now study which inequality do we get with the word length. The cocycle associated to it is infinite-dimensional, with an orthonormal basis which can be described as oriented edges in the coordinate axes of the Cayley graph of $\Z^n$. More precisely, the associated Gromov form on $\R[\Z^n]$ is
\begin{align*}
	\langle \delta_g, \delta_h \rangle_{\psi_1} = \frac{1}{2} \big( \psi_1(g) + \psi_1(h) - \psi_1(h-g) \big) = \Suma{j=1}{n} \min \big\{ |h_j|,|g_j| \big\} \, \delta_{g_j \cdot h_j > 0}.
\end{align*}
Given $g \in \Z^n$ and $j \in [n]$, define $$g_{[j]}^- = g - \mathrm{sgn}(g_j) e_j \quad \mbox{with} \quad \mathrm{sgn}(0)=0.$$ Then, we may construct the following elements in $\R[\Z^n]$ $$w_{g,j} = \delta_g - \delta_{g_{[j]}^-} \quad \mbox{and} \quad u_j(\ell) = w_{\ell e_j,j}.$$
If $\mathcal{H}_{\psi_1} = \R[\Z^n]/\mathrm{Ker} \langle \cdot, \cdot \rangle_{\psi_1}$, the following properties define an orthonormal basis: 
\begin{itemize}
\item $\langle u_j(\ell), u_j(\ell) \rangle_{\psi_1} = 1$ for all $(j,\ell) \in [n] \times \Z \! \setminus \! \{0\}$.

\vskip3pt

\item $\langle u_j(\ell), u_{j'}(\ell') \rangle_{\psi_1} = 0$ whenever $j \neq j'$ or $\ell \neq \ell'$.

\vskip3pt

\item $w_{g,j} = u_j(g_j)$ since the difference belongs to $\mathrm{Ker} \langle \cdot, \cdot \rangle_{\psi_1}$.
\end{itemize}
Altogether, this implies that the set
$$
\Big\{ u_j(\ell): (j,\ell) \in [n] \times \Z \! \setminus \! \{0\} \Big\}
$$
is an orthonormal basis for $\mathcal{H}_{\psi_1}$. The cocycle map is given by $\beta(g) = \delta_{g}$ and the orthogonal action $\alpha$ satisfies $\alpha_g(\delta_h) = \delta_{g+h} - \delta_{g}$. This means that for any $g \in \Z^n$ we have 
$$
\alpha_g(u_j(\ell)) = \delta_{g + \ell e_j} - \delta_{g + \ell e_j - (\mathrm{sgn}(\ell))e_j}.
$$
Therefore, the subspaces $\mathcal{H}_{\psi_1,j} = \mbox{span} \{u_j(\ell): \ell \in \Z \setminus \{0\}\}$ are $\alpha$-invariant 
for $j \in [n]$. This proves that the same conditional expectations $\cexp_{[n]\setminus \Ss}$ considered before still define an admissible family of Fourier truncations. The cocycle derivative associated to $u_j(\ell)$ acts as follows:
\begin{align*}
	\partial_{u_j(\ell)} \lambda(g) & = 2\pi i \, \langle u_j(\ell), \delta_g \rangle_{\psi_1} \lambda(g) \\
	&= 2\pi i \, \Big( \min \big\{ |g_j|,|\ell| \big\} \delta_{g_j \cdot \ell > 0} - \min \big\{ |g_j|,|\ell|-1 \big\} \delta_{g_j \cdot (\ell - \mathrm{sgn}(\ell)) > 0} \Big) \lambda(g) \\
	&= 2\pi i \, \delta_{\{g_j \cdot \ell > 0, |g_j| \geq |\ell|\}} \lambda(g).
\end{align*}
The Laplacian is 
\begin{equation*}
\laplacian_{\psi_1} f = \suma{g \in \Z^n} \psi_1(g) \widehat{f}(g) \lambda(g),
\end{equation*}
whose spectral gap is still $\sigma = \min_j \psi_1(e_j)=1$. Theorem \ref{thm:theoremA} yields 
\begin{equation}\label{eqn:thAZ1}
\frac{1}{{n \choose k}} \suma{\substack{\Ss \subseteq [n] \\ |\Ss|=k}} \Big\| \suma{g \in \Z^\Ss}^{\null} \widehat{f}(g) e^{2\pi i \langle \cdot, g \rangle} \Big\|^p_{L_p(\mathbb{T}^n)} \lesssim_p \frac{k}{n} \Suma{j = 1}{n} \| \gradient_j f \|_{L_p(\mathbb{T}^n)}^p + \Big(\frac{k}{n}\Big)^{\frac{p}{2}} \|f\|^p_{L_p(\mathbb{T}^n)},
\end{equation}
with 
$$
\| \gradient_j (f) \|_{L_p(\mathbb{T}^n)} = \| \gradient_j (f^*) \|_{L_p(\mathbb{T}^n)} = \Big\| \Big( \sum_{\ell \in \mathbb{Z} \setminus \! \{0\}} |\partial_{u_j(\ell)} f|^2\Big)^{\frac{1}{2}} \Big\|_{L_p(\mathbb{T}^n)}.
$$

\begin{remark}
\emph{Note that $|\partial_{u_j(\ell)}(f)| \neq |\partial_{u_j(\ell)}(f^*)|$. Thus, nontrivial cocycle actions lead to noncommutative phenomena even when working with abelian groups, as pointed out in \cite{JMP18}. In spite of that, observe that $\langle \delta_{-g}, u_j(\ell) \rangle_\psi = \langle \delta_g, u_j(-\ell) \rangle_\psi$ which implies $\|\gradient_j(f)\|_p = \|\gradient_j(f^*)\|_p$ as claimed above.}
\end{remark}

\noindent On the other hand, taking 
\begin{equation*}
	\partial_j \lambda(g) := \frac{1}{2\pi i} \Big(\partial_{u_j(1)} + \partial_{u_j(-1)} \Big) = \delta_{g_j \neq 0} \lambda(g)
\end{equation*}
we get $\partial_{u_j(\ell)} \circ \partial_j = \partial_{u_j(\ell)}$ for any $(j,\ell) \in [n] \times \Z \! \setminus \! \{0\}$. Thus, Theorem \ref{thm:theoremC} gives
\begin{equation*}
	\frac{1}{{n \choose k}} \suma{\substack{\Ss \subseteq [n] \\ |\Ss|=k}} \Big\| \suma{g \in \Z^\Ss}^{\null} \widehat{f}(g) e^{2\pi i \langle \cdot, g \rangle} \Big\|^p_{L_p(\mathbb{T}^n)} \lesssim_p \frac{k}{n} \Suma{j=1}{n} \|\partial_j f\|_{L_p(\mathbb{T}^n)}^p + \Big( \frac{k}{n} \Big)^{\frac{p}{2}} \|f\|_{L_p(\mathbb{T}^n)}^p
\end{equation*}
for any mean-zero $f \in L_p(\mathbb{T}^n)$. This recovers inequality \eqref{eq:ZntheoremC} and improves \eqref{eqn:thAZ1}.


\subsection{Discrete tori} 

Consider the word length $|g|= \min\{g,2m-g\}$ in $\mathbb{Z}_{2m}$. As shown in \cite{JPPP17}, it defines a conditionally negative symmetric length. In particular the same holds for the corresponding length in the product $\mathbb{Z}_{2m}^n$ 
\begin{equation*}
	\psi(g) = |g_1| + |g_2| + \ldots + |g_n| \quad \mbox{for} \quad g=(g_1,\ldots,g_n) \in \Z_{2m}^n.	
\end{equation*}
This word length has many similarities with the previous one
\begin{equation*}
\langle \delta_g, \delta_h \rangle_\psi = \frac{1}{2} \big( \psi(g)+\psi(h)-\psi(h-g) \big) = \frac{1}{2} \Suma{j=1}{n} |g_j| + |h_j| - |h_j - g_j|.
\end{equation*}
Given $g \in \Z_{2m}^n$ and $j \in [n]$, define $$w_{g,j} = \delta_g - \delta_{g - e_j} \quad \mbox{and} \quad u_j(\ell) = w_{\ell e_j,j} \quad \mbox{for} \quad 1 \le \ell \le 2m.$$
If $\mathcal{H}_{\psi} = \R[\Z_{2m}^n]/\mathrm{Ker} \langle \cdot, \cdot \rangle_{\psi}$, we find that 
\begin{itemize}
\item $\langle u_j(\ell), u_j(\ell) \rangle_{\psi} = 1$ for all $(j,\ell) \in [n] \times [m]$.

\vskip3pt

\item $\langle u_j(\ell), u_{j'}(\ell') \rangle_{\psi} = 0$ whenever $j \neq j'$ or $\ell \neq \ell', \ell'+m$.

\vskip3pt

\item $w_{g,j} = u_j(g_j)$ since the difference belongs to $\mathrm{Ker} \langle \cdot, \cdot \rangle_{\psi}$.

\vskip3pt

\item $u_j(\ell) = - u_j(\ell+m)$ since the difference belongs to $\mathrm{Ker} \langle \cdot, \cdot \rangle_{\psi}$.

\vskip3pt

\item $\langle \delta_{\ell e_j}, \delta_{\ell' e_j} \rangle_\psi = \min \big\{ \ell, 2m-\ell', \max\{0,m-\ell'+\ell\} \big\}$ for $1\le \ell \leq \ell' \le 2m$.
\end{itemize}
Altogether, this implies that the set
$$
\Big\{ u_j(\ell): (j,\ell) \in [n] \times [m] \Big\}
$$
is an orthonormal basis for $\mathcal{H}_{\psi}$. The cocycle map is given by $\beta(g) = \delta_{g}$ and the orthogonal action $\alpha$ satisfies $\alpha_g(\delta_h) = \delta_{g+h} - \delta_{g}$. This means that for any $g \in \Z^n_{2m}$ we have 
$$
\alpha_g(u_j(\ell)) = \delta_{g + \ell e_j} - \delta_{g + (\ell - 1)e_j}.
$$
Therefore, the subspaces $\mathcal{H}_{\psi,j} = \mbox{span}\{u_j(\ell): \ell \in [m]\}$ give again an $\alpha$-invariant  splitting of $\mathcal{H}_\psi$ with $j$ running over $[n]$. In particular, the conditional expectations $\cexp_{[n]\setminus \Ss}$ over the subgroups $\Z_{2m}^\mathsf{S}$ define an admissible family of Fourier truncations and the cocycle derivatives are given by
$$
\partial_{u_j(\ell)} \lambda(g) = 2\pi i \ \delta_{\{\ell \leq g_j < \ell+m\}} \lambda(g).
$$
The associated Laplacian has spectral gap equal to $1$. As before, Theorem \ref{thm:theoremA} yields a statement that we omit because it can readily be improved. If we set as before $\partial_j \lambda(g):= \delta_{g_j \neq 0} \lambda(g) $ for $j \in [n]$, we immediately see that $\partial_{u_j(\ell)} \circ \partial_j = \partial_{u_j(\ell)}$. Moreover, we can rewrite it as follows 
$$
\partial_j \lambda(g) =  \frac{1}{2\pi i} \Big( \partial_{u_j(1)} \lambda(g) + \partial_{u_j(m)} \lambda(g) \Big) - \delta_{g_j = m}  \lambda(g).
$$
Next, note that we may write the last term as follows
$$
\delta_{g_j = m} \lambda(g) = \cexp_{\{0,m\},j} \Big(\frac{1}{2\pi i}\partial_{u_j(1)} \lambda(g) \Big) = \frac{1}{4\pi i} \cexp_{\{0,m\},j} \Big( \partial_{u_j(1)} \lambda(g) + \partial_{u_j(m)} \lambda(g) \Big),
$$
where $\cexp_{\{0,m\},j}$ is the conditional expectation onto $\Z_{2m}^{j-1} \times \{0,m\}\times \Z_{2m}^{n-j}$. Then, after applying Theorem \ref{thm:theoremC} one gets the following for mean-zero $f: \Z_{2m}^n \to \C$ 
\begin{eqnarray} \label{eq:thmBforZ_2m^n}
\hskip20pt \frac{1}{{n \choose k}} \suma{\substack{\Ss \subseteq [n] \\ |\Ss|=k}}  \Big\| \sum_{g \in \Z_{2m}^\mathsf{S}} \widehat{f}(g) e^{\frac{\pi i}{m} \langle \cdot, g \rangle} \Big\|_p^p \!\!\!\! & \lesssim_p & \!\!\!\! \frac{k}{n} \Suma{j=1}{n} \|\partial_j f\|_p^p + \Big( \frac{k}{n} \Big)^{\frac{p}{2}} \|f\|_p^p \\ \nonumber \!\!\!\! & \lesssim & \!\!\!\! \frac{k}{n} \Suma{j=1}{n} \big\| ( \partial_{u_j(1)} + \partial_{u_j(m)} )f \big\|_p^p + \Big( \frac{k}{n} \Big)^{\frac{p}{2}} \|f\|_p^p.
\end{eqnarray}

\begin{remark}
\emph{Inequality \eqref{eq:thmBforZ_2m^n} for $\Z_{2m}^n$ with $m=1$ generalizes Naor's inequality \eqref{Eq-Naor_p} for the hypercube. Just identify $g \in \Z_2$ with $\exp(\pi i g) \in \{\pm 1\}$. Then note that $\partial_j^1 = 2 \partial_j^2$, where $\partial_j^1$ is the discrete derivative used by Naor and $\partial_j^2$ is our choice of $\partial_j$ in \eqref{eq:thmBforZ_2m^n} for $m=1$. Moreover, we can consider weighted forms of Naor's inequality by considering different measures on the same group $\mathbb{Z}_2^n$ to get different cocycle representations. We next show that this does not lead to an improvement over the result in \cite{N16}. Indeed, let $\G=\{-1,1\}^n$ and equip  $\Gamma = \widehat{\G} = \{-1,1\}^n$ with the measure
\begin{equation*}
\mu = \Suma{j=1}{n} \alpha_j \delta_{w_j},
\end{equation*}
with $\alpha_j \geq 0$ and $w_j = (1, \ldots, 1, -1, 1, \ldots, 1) = (1, 1, \ldots, 1) - 2e_j$ for $j \in [n]$. Identifying $\Gamma$ with the power set of $[n]$ we identify $w_j$ with $\{j\}$. Following Example B from \cite[Subsection 1.4]{JMP18}, we consider the conditionally negative length function
\begin{equation*}
	\psi(\mathsf{A}) := \big\| 1 - W_\mathsf{A} \big\|^2_{L_2(\Gamma,\mu)}.
\end{equation*}
Then $\psi$ may be represented by the cocycle $(\mathcal{H}_\psi, \alpha, \beta)$ with $$\mathcal{H}_\psi = L_2(\Gamma,\mu), \quad \alpha_\mathsf{A}(u) = W_\mathsf{A} \cdot u, \quad \beta(\mathsf{A})= 1 - W_\mathsf{A}.$$ Then $\big\{ u_j = \alpha_j^{-\frac12}  \delta_{\{j\}}: j \in [n] \big\}$ is an ONB and Riesz transforms take the form 
\begin{equation*}
R_{u_j}f = \suma{\mathsf{A} \subset [n]} \frac{\langle \beta(\mathsf{A}),\delta_{\{j\}}\rangle_\psi}{\sqrt{\alpha_j \psi(\mathsf{A})}}  \widehat{f}(\mathsf{A}) W_\mathsf{A} = \suma{\substack{\mathsf{A} \subset [n] \\ j \in \mathsf{A}}} \frac{\sqrt{\alpha_j}}{\sqrt{\sum_{\ell \in A} \alpha_\ell}} \widehat{f}(\mathsf{A}) W_\mathsf{A}.
\end{equation*}
Consider the decomposition $\mathcal{H}_{\psi,j} = \mathbb{R} \delta_{\{j\}}$. Note $\alpha_\mathsf{A}(\delta_{\{j\}}) = W_\mathsf{A} \delta_{\{j\}} = (-1)^{\delta_{j \in \mathsf{A}}} \delta_{\{j\}}$, so $\alpha_\mathsf{A}(\mathcal{H}_{\psi,j}) \subset \mathcal{H}_{\psi,j}$ and the decomposition is equivariant. Therefore, the associated conditional expectation can be chosen to be
\begin{equation*}
	\cexp_{[n]\setminus \Ss} f = \suma{\beta(\mathsf{A}) \in \mathcal{H}_\Ss} \widehat{f}(\mathsf{A}) W_\mathsf{A} \, = \, \suma{\mathsf{A} \subset \Ss} \widehat{f}(\mathsf{A}) W_\mathsf{A}.
\end{equation*}
Finally, the cocycle derivatives are given by 
\begin{eqnarray*}
\partial_{u_j} W_\mathsf{A} \!\!\! & = & \!\!\! \frac{2\pi i}{\sqrt{\alpha_j}}  \langle \beta(\mathsf{A}), \delta_{\{j\}} \rangle_\psi W_\mathsf{A} \\ \!\!\! & = & \!\!\! 4 \pi i \sqrt{\alpha_j} \delta_{j \in \mathsf{A}} W_\mathsf{A} = 2\pi i \sqrt{\alpha_j} \partial_j W_\mathsf{A},
\end{eqnarray*} 
where $\partial_j$ denotes the $j$-th discrete derivative. Therefore, Theorem \ref{thm:theoremA} yields 
\begin{eqnarray*}
\lefteqn{\hskip-40pt \frac{1}{{n \choose k}} \suma{\substack{\Ss \subseteq [n] \\ |\Ss|=k}} \Big\| \suma{\mathsf{A} \subset \Ss} \widehat{f}(\mathsf{A}) W_\mathsf{A} \Big\|^p_{L_p(\Omega_n)}} \\ \hskip40pt \!\!\! & \lesssim & \!\!\! \frac{1}{\sigma^{p/2}} \frac{k}{n} \Suma{j=1}{n} \alpha_j^{\frac{p}{2}} \|\partial_j f\|^p_{L_p(\Omega_n)} + \Big(\frac{k}{n}\Big)^{\frac{p}{2}} \|f\|^p_{L_p(\Omega_n)} \\ \!\!\! & = & \!\!\! \frac{k}{n} \Suma{j=1}{n} \Big( \frac{\alpha_j}{\min_{k\in [n]} \alpha_k}\Big)^{\frac{p}{2}} \|\partial_j f\|^p_{L_p(\Omega_n)} + \Big( \frac{k}{n} \Big)^{\frac{p}{2}} \|f\|_{L_p(\Omega_n)}^p,
\end{eqnarray*}
since $\sigma = \min_{k \in [n]} \psi(\{k\}) = 4 \min_{k \in [n]} \alpha_k$. Thus taking $\alpha_j = 1$ for all $j$ is optimal.}
\end{remark}

\begin{remark}
\emph{It is natural to ask if the situation changes much when the cyclic groups under consideration have odd cardinal. The function $\psi(g) = \sum_{j=1}^n |g_j|$ with $|g_j| = \min\{g_j, 2m+1-g_j\}$ is a conditionally negative length on $\Z_{2m+1}^n$, and so there exists an associated cocycle induced by the Gromov form
\begin{align*}
\langle \delta_g, \delta_h \rangle & = \frac{1}{2} \big( \psi(g) + \psi(h) - \psi(g-h) \big) \\
& = \sum_{j=1}^n \min \Big\{ g_j,2m+1-h_j,\max \big\{0,m-h_j+g_j+\frac{1}{2} \big\} \Big\}.
\end{align*}
It defines a cocycle Hilbert space $\mathcal{H}_\psi$ with dimension $2mn$. Theorems \ref{thm:theoremA} and \ref{thm:theoremC} apply but calculating an explicit expression for the orthonormal basis of $\mathcal{H}_\psi$ is more tedious and we shall leave it to the interested reader.}
\end{remark}


\section{\bf Applications to free products} \label{sec3} 

We now explore applications of Theorem \ref{thm:theoremC} after replacing the direct products in the previous section by free products. Given a free product $\G = \G_1 * \G_2 * \cdots * \G_n$ a general element $g \in \G$ can always be written in reduced form $g = g_{i_1} g_{i_2} \cdots g_{i_s}$ where $g_{i_k} \in \G_{i_k}$ and $i_1 \neq i_2 \neq \cdots \neq i_s$. We shall be working with the free group $\F_n = \Z * \Z * \cdots * \Z$ and with the free product $\Z_{2m}^{*n}$ of $n$ copies of $\Z_{2m}$. In both cases we shall write $g_1, g_2, \ldots, g_n$ for the canonical generators and a generic element will be a word of the form $$w = g_{i_1}^{\ell_1} g_{i_2}^{\ell_2} \cdots g_{i_s}^{\ell_s}$$ with $i_1 \neq i_2 \neq \cdots \neq i_s$ and $\ell_k$ in $\Z$ or $\Z_{2m}$ accordingly.  

\subsection{The free group} 

Define 
$$
|w| = \Suma{j=1}{r} |\ell_j| \quad \mbox{for} \quad w = g_{i_1}^{\ell_1} \ldots g_{i_r}^{\ell_r} .
$$
Haagerup proved in \cite{H78} that it is conditionally negative. The cocycle structure naturally induced by the word length $|\cdot|$ can be described through the Hilbert space orthonormaly generated by outgoing oriented edges in its Cayley graph. To be more precise, let us consider the following partial order on $\mathbb{F}_n$. Given $w_1 = g_{i_1}^{\ell_1} \ldots g_{i_r}^{\ell_r}$ and $w_2 = g_{j_1}^{t_1} \ldots g_{j_s}^{t_s}$ with $\ell_j, t_j \in \Z \setminus \{0\}$, we say that $w_1 \leq w_2$ when 
\begin{itemize}
	\item $r \leq s$,
	\item $g_{i_k}^{\ell_k} = g_{j_k}^{t_k} \ \mbox{for } 1 \leq k \leq r-1$, 
	\item $g_{i_r} = g_{j_r}$, $\ell_{k} t_k > 0$ and $|\ell_k| \leq |t_k|$.
\end{itemize}
Any $w_1 \leq w_2$ is called an initial subchain of $w_2$. As we did with elements of cyclic groups equipped with their natural order structure, we can now define predecessors. If $w = g_{i_1}^{\ell_1} \ldots g_{i_r}^{\ell_r} \neq e$, we denote 
\begin{equation*}
	w^{-} = g_{i_1}^{\ell_1} \ldots g_{i_{r}}^{\ell_r - \mathrm{sgn}(\ell_r)}.
\end{equation*}
The Gromov form takes the following form in this case
\begin{equation*}
\langle \delta_{w_1} , \delta_{w_2}  \rangle_{|\cdot|} = \frac{1}{2} \big( |w_1| + |w_2| - |w_1^{-1}w_2| \big) = |\min \{w_1,w_2\}|,
\end{equation*}
where $\min \{w_1,w_2\}$ denotes the longest word which is an initial chain of both $w_1$ and $w_2$. Given $w \neq e$ in $\F_n$, we define $u_w = \delta_w - \delta_{w^-} \in \R[\F_n]$. Very much like in the previous section, we find the following properties:
\begin{itemize}
\item $\mathrm{Ker}(	\langle \cdot , \cdot \rangle_{|\cdot|}) = \R \delta_e$.

\vskip3pt

\item $\langle u_w , u_w \rangle_{|\cdot|} = 1$ for $w \in \F_n \setminus \{e\}$. 

\vskip3pt

\item $\langle u_{w_1} , u_{w_2} \rangle_{|\cdot|} = 0$ for $w_1 \neq w_2$ in $\F_n$. 
\end{itemize}
This proves that $$\Big\{ u_w: w \in \F_n \setminus \{e\} \Big\}$$ is an orthonormal basis of $\mathcal{H}_{|\cdot|} = \R[\F_n] / \R\delta_e$. The cocycle map and the cocycle action are determined as usual by $\beta(w) = \delta_w$ and $\alpha_w(\delta_{w'}) = \delta_{ww'} - \delta_w$.The cocycle derivative in the direction of $u_{w}$ is
\begin{equation*}
	\partial_{u_w} \lambda(w') = 2\pi i \langle \beta(w'), u_{w} \rangle \lambda(w') = 2\pi i \delta_{w \leq w'} \lambda(w') \ \Rightarrow \ \partial_{u_w} f = 2\pi i \hskip-3pt \suma{w \leq w'} \widehat{f}(w') \lambda(w').
\end{equation*}
Next, we decompose $\mathcal{H}_{|\cdot|}$ as
$$
\mathcal{H}_{|\cdot|} = \bigoplus_{j=1}^n \mathcal{H}_{|\cdot|,j} \quad \mbox{with} \quad \mathcal{H}_{|\cdot|,j}= \mbox{span} \big\{u_w: g_j \le w \mbox{ or } g_j^{-1} \leq w \big\}.
$$ 
This leads to consider the Fourier truncations 
$$
\cexp_{[n]\setminus \Ss} f := \suma{w \in \mathbb{F}_\Ss}^{\null} \widehat{f}(w) \lambda(w).
$$
Being conditional expectations, these Fourier truncations are completely contractive and pairwise $\beta$-orthogonality holds since we trivially have $\beta(\mathbb{F}_\Ss) = \beta(\mathbb{F}_\Ss^{-1}) \subset \mathcal{H}_{|\cdot|,\Ss}$. Next, taking the derivatives 
$$
\partial_j :=  \frac{1}{2\pi i} \Big(\partial_{u_{g_j}} + \partial_{u_{g_j^{-1}}}\Big) \quad \mbox{for} \quad j \in [n],
$$
we can readily check that $\partial_{u_w} \circ \partial_{j_1}  = \partial_{u_w}$ whenever $u_w \in \mathcal{H}_{|\cdot|,j_1}$. In conclusion, we have checked all the hyporthesis to apply Theorem \ref{thm:theoremC} for our family of Fourier truncations. In this case we get
\begin{equation} \label{eq:free1}
\frac{1}{{{n}\choose{k}}} \sum_{\begin{subarray}{c} \mathsf{S} \hskip1pt \subset \hskip1pt [n] \\ |\mathsf{S}|=k \end{subarray}} \Big\| \sum_{w \in \F_\mathsf{S}} \widehat{f}(w) \lambda(w) \Big\|_p^p \, \lesssim_p \, \frac{k}{n} \sum_{j=1}^n \Big[\big\| \partial_j(f) \big\|_p^p + \big\| \partial_j(f^*) \big\|_p^p \Big]+ \Big( \frac{k}{n} \Big)^{\frac{p}{2}} \|f\|_p^p.
\end{equation}
Inequality \eqref{eq:free1} is very close to the conjectured free form of Naor's inequality \eqref{Eq-FreeNaor_p} in the Introduction, with an extra adjoint  term which we shall eliminate at the end of the paper by proving an even stronger inequality. 



\subsection{The free product $\Z_{2m}^{*n}$} A similar analysis applies as well in this case. Given two reduced words $w_1 = g_{i_1}^{\ell_1} \ldots g_{i_r}^{\ell_r}$ and $w_2 = g_{j_1}^{t_1} \ldots g_{j_s}^{t_s}$ with $\ell_j, t_j \in [2m-1]$, we say that $w_1 \leq w_2$ if and only if
\begin{itemize}
	\item $r \leq s$,
	\item $i_k = j_k$ for any $k \in [r]$ and $\ell_k = t_k$ for any $k \in [r-1]$,
	\item either $\ell_r, t_r \in [m]$ and $i_r \leq j_r$, or $i_r, j_r \in [2m-1] \setminus [m-1]$ and $i_r \geq j_r$.
\end{itemize}
The map $\psi:\Z_{2m}^{*n} \to \mathbb{R}_+$ given by 
\begin{equation*}
	w= g_{i_1}^{\ell_1} \ldots g_{i_r}^{\ell_r} \mapsto \psi(w)=\sum_{k=1}^r |g_{i_j}^{\ell_j}|=\Suma{k=1}{r} \min\{\ell_k,2m-\ell_k\}
\end{equation*}
is a conditionally negative length function \cite{H78}, with associated Gromov form
\begin{eqnarray}\label{eq:FreeZ2mnGromov}
\langle \delta_{w_1}, \delta_{w_2} \rangle_{\psi} \!\!\! & = & \!\!\! \frac{1}{2} \big( \psi(w_1) + \psi(w_2) - \psi(w_1^{-1}w_2) \big) \\ \nonumber \!\!\! & = & \!\!\! \psi \big( \min \{w_1,w_2\} \big) + \frac{1}{2} \big( \psi(\eta_1)+\psi(\eta_2)-\psi(\eta_1^{-1} \eta_2) \big).
\end{eqnarray}
where $\min \{w_1,w_2\}$ is again the longest common subchain and $w_j = \min \{w_1, w_2\} \eta_j$ for $j=1,2$. The second term above is always $0$ in the free group $\F_n$, but not necessarily in this case. Given $w= g_{i_1}^{\ell_1} \ldots g_{i_r}^{\ell_r} \neq e$ we define $w^-=g_{i_1}^{\ell_{1}} \ldots g_{i_r}^{\ell_r -1}$ and construct $u_w = \delta_w - \delta_{w^-}$ as usual. Then, we find that
\begin{itemize}
\item $\langle u_{w}, u_{w} \rangle_{\psi} = 1$ for every $w \in \Z_{2m}^{*n} \setminus \{e\}$.

\vskip1pt

\item $\langle u_{w_1}, u_{w_2} \rangle_{\psi} = 0$ when $e \neq w_1^{-1}w_2 \neq g_j^m$ for $j \in [n]$.

\item $\langle u_{w_1}, u_{w_2} \rangle_{\psi} = 0$ when $w_1^{-1}w_2 = g_j^m$ and both $w_1,w_2$ end with $g_j^{\pm1}$.


\item If $w = g_{i_1}^{\ell_1} \ldots g_{i_r}^{\ell_r}$, then $u_w = - u_{wg_{i_r}^m}$ in $\mathcal{H}_{\psi} = \R[\Z_{2m}^{*n}]/ \mathrm{Ker} \langle \cdot, \cdot \rangle_{\psi}$.
\end{itemize}
This proves that $$\Big\{ u_w: w = g_{i_1}^{\ell_1} \ldots g_{i_r}^{\ell_r} \in \Z_{2m}^{*n} \setminus \{e\}  \ \mathrm{ with } \ \ell_r \in [m] \Big\}$$ is an orthonormal basis of $\mathcal{H}_{\psi} = \R[\Z_{2m}^{*n}] / \mathrm{Ker} \langle \cdot, \cdot \rangle_{\psi}$. We set as usual $\beta(w) = \delta_w$ and $\alpha_w(\delta_{w'}) = \delta_{ww'} - \delta_w$. Among the above properties it is perhaps convenient to justify the last one. Note that $\langle u_w + u_{wg_{i_r}^m}, u_w + u_{wg_{i_r}^m} \rangle_{\psi} = 0$ iff  $\langle u_w, u_{wg_{i_r}^m} \rangle_{\psi} = -1$ but we have
\begin{eqnarray*}
\langle u_w, u_{wg_{i_r}^m} \rangle_{\psi} \!\!\! & = & \!\!\! \frac12 \Big( -\psi \big(g_{i_r}^m\big) + \psi\big((w^{-})^{-1}wg_{i_r}^m\big) + \psi\big(g_{i_r}^{m-1}\big) - \psi\big((w^-)^{-1} w g_{i_r}^{m-1}\big) \Big) \\ \!\!\! & = & \!\!\! \frac12 \Big( -\psi(g_{i_r}^m) + \psi(g_{i_r}^{m+1}) + \psi(g_{i_r}^{m-1}) - \psi(g_{i_r}^{m}) \Big) \\ \!\!\! & = & \!\!\! \frac12 \Big( -m + m-1 + m-1 -m \Big) \, = \, -1.
\end{eqnarray*}
If $w = g_{i_1}^{\ell_1} \ldots g_{i_r}^{\ell_r}$ with $\ell_r \in [m]$, derivatives are given by 
\begin{equation} \label{eq:DerivativeW}
\partial_{u_w} \lambda(w^\prime) = 2\pi i \, \langle \beta(w^\prime) , u_w \rangle_\psi \lambda(w^\prime) = 2\pi i \, \delta_{w^\prime \in W(w)} \lambda(w^\prime)
\end{equation}
where $W(w)$ is the set of those words $w' = g_{j_1}^{t_1} \ldots g_{j_s}^{t_s}$ satisfying
\begin{equation} \label{eq:SetW}
r \leq s, \quad i_k = j_k \ \mathrm{for \ } k\le r, \quad \ell_k = t_k \ \mathrm{for \ } k \le r-1 \quad \mbox{and} \quad \ell_r \leq t_r \leq \ell_r + m -1.
\end{equation}
Indeed, just write $\beta(w') = \delta_{w'} = u_{w'} + \delta_{{w'}^{-}} = u_{w'} + u_{{w'}^{-}} + \delta_{{w'}^{--}}$ and so on. The inner product with $u_w$ will be $0$ unless we find $u_w$ in our telescopic sum above just once, in which case we get the value $1$. Note that it could appear twice due to the identity $u_w = - u_{wg_{i_r}^m}$ recalled above. In that case, they get mutually cancelled and we get $0$. This happens when $t_r-\ell_r \in [2m-1] \setminus [m-1]$. 

It remains to consider Fourier truncations. As for the free group, our choice is the conditional expectation into the subgroup $\mathbb{Z}_{2m}^{*\Ss} = \langle g_{j} : j \in \Ss \rangle$. Then we consider the decomposition $$\mathcal{H}_\psi = \bigoplus_{j=1}^n \mathcal{H}_{\psi,j} \quad \mbox{with} \quad \mathcal{H}_{\psi,j} = \mathrm{span} \big\{ u_w: g_j \le w \ \mbox{or} \ g_j^{-1} \le w \big\}.$$
Our Fourier truncations form an admissible family. Define 
$$
\partial_j \lambda(w) = \frac{1}{2\pi i} \Big( \partial_{u_{g_j}} \lambda(w) + \partial_{u_{g_j^m}} \lambda(w) \Big) - \delta_{g_{i_1}^{\ell_1} = g_{j}^m} \lambda(w) \quad \mbox{for any} \quad w = g_{i_1}^{\ell_1} \ldots g_{i_r}^{\ell_r}.
$$
In other words, $\partial_j \lambda(w) = \delta_{i_1=j} \lambda(w)$ for $w\neq e$ and $\partial_{u_w} \circ \partial_{j} = \partial_{u_w}$ for $u_w \in \mathcal{H}_{\psi,j}$. The construction above yields the form of Theorem \ref{thm:theoremC} on the von Neumann algebra of the free product $\Z_{2m}^{*n}$
$$
\frac{1}{{n \choose k}} \suma{\substack{\Ss \subseteq [n] \\ |\Ss|=k}}  \Big\| \sum_{w \in \Z_{2m}^{*\Ss}} \widehat{f}(w) \lambda(w) \Big\|_p^p \lesssim_p \frac{k}{n} \Suma{j=1}{n} \big\| \partial_j (f) \big\|_p^p + \big\| \partial_j (f^*) \big\|_p^p + \Big( \frac{k}{n} \Big)^{\frac{p}{2}} \|f\|_p^p.
$$


\subsection{Free Hilbert transforms} \label{subsec:freeHilbert} 

Compared to \eqref{Eq-FreeNaor_p}, the form of Theorem \ref{thm:theoremC} for free groups gives the additional terms $\partial_j(f^*)$. These terms seem to be necessary in the general context of Theorem \ref{thm:theoremC}, but they are removable for free groups |in fact, we shall prove an even stronger inequality| due to a singular behavior of word-length derivatives for free groups. This comes from the following identity $$\langle \delta_{w_1}, u_{w_2} \rangle_{|\cdot|} = \mathrm{sgn} \big( \langle \delta_{w_1}, u_{w_2} \rangle_{|\cdot|} \big),$$ since the above inner product can only take the values $0$ or $1$. This means in particular that word-length derivatives can be regarded as free forms of directional Hilbert transforms, which were recently investigated by Mei and Ricard in \cite{MR17}. To be more precise, if $\mathbb{A}_\Ss$ is the set of free words whose first letter is in $\F_\Ss$, it turns out that distinguished derivatives satisfy $$\partial_j = \frac{1}{2\pi i} \Big( \partial_{u_{g_j}} + \partial_{u_{g_j^{-1}}} \Big) = \mbox{Projection onto $\mathbb{A}_{\{j\}}$}.$$ The free Hilbert transforms for mean-zero $f$ are defined as 
$$
H_{\varepsilon}(f) = \Suma{j = 1}{n} \varepsilon_j \partial_j(f) \quad \mbox{for} \quad \varepsilon_j = \pm 1.
$$
Mei and Ricard proved in \cite{MR17} the following crucial inequality 
\begin{equation}\label{eq:FreeHilbertTransform}
\|H_\varepsilon f\|_{L_p(\mathcal{L}(\mathbb{F}_n))} \asymp_p \|f\|_{L_p(\mathcal{L}(\mathbb{F}_n))} \quad \mbox{for any} \quad 1 < p < \infty.
\end{equation}

\begin{theorem} \label{FreeXpChaosHilbert}
If $p \ge 2$ and $k \in [n]$, every mean-zero $f \in L_p(\mathcal{L}(\mathbb{F}_n))$ satisfies
\begin{equation*}
\frac{1}{{n \choose k}} \suma{\substack{\Ss \subseteq [n] \\ |\Ss|=k}} \Big\| \sum_{w \in \mathbb{A}_\Ss} \widehat{f}(w) \lambda(w) \Big\|_{L_p(\mathcal{L}(\mathbb{F}_n))}^p \lesssim_p \frac{k}{n} \Suma{j=1}{n} \big\| \partial_j (f) \big\|_{L_p(\mathcal{L}(\mathbb{F}_n))}^p + \Big(\frac{k}{n}\Big)^{\frac{p}{2}} \|f\|_{L_p(\mathcal{L}(\mathbb{F}_n))}^p.
\end{equation*}
\end{theorem}

\dem	Define $$h = \sum_{w \in \mathbb{A}_\Ss} \widehat{f}(w) \lambda(w) = \sum_{j \in \Ss} \partial_j (f).$$
Applying inequality \eqref{eq:FreeHilbertTransform} we obtain
$$
\|h\|_p \asymp_p \mathbb{E} \big\| H_\varepsilon(h) \big\|_p = \mathbb{E} \Big\| \sum_{j \in \Ss} \varepsilon_j \partial_j(f) \Big\|_p.
$$
The result follows from Theorem \ref{thm:linearXp} and another application of \eqref{eq:FreeHilbertTransform} for $f$. \fin 

\begin{corollary} \label{Cor-FNq}
Inequality \eqref{Eq-FreeNaor_p} holds for $p \ge 2$ and any mean-zero $f \in L_p(\mathcal{L}(\F_n))$.
\end{corollary}

\dem It follows from Theorem \ref{FreeXpChaosHilbert} and the trivial inequality 
$$
\Big\| \sum_{w \in \F_\Ss} \widehat{f}(w) \lambda(w) \Big\|_p \, = \, \Big\| \sum_{w \in \F_\Ss} \widehat{h}(w) \lambda(w) \Big\|_p \, \le \, \|h\|_p \, = \, \Big\| \sum_{w \in \mathbb{A}_\Ss} \widehat{f}(w) \lambda(w) \Big\|_p,
$$
where $h$ is defined as in the proof of Theorem \ref{FreeXpChaosHilbert}, since we note that $\F_\Ss \subset \mathbb{A}_\Ss$. \fin

\begin{remark}
\emph{It is conceivable that Theorem \ref{FreeXpChaosHilbert} or at least Corollary \ref{Cor-FNq} could have been proved as well from a generalized form of Theorem \ref{thm:theoremC} in the line of Remark \ref{Rem-ThmB'}, but we have not found an argument using such an approach.}  
\end{remark}

\begin{remark}
\emph{Hilbert transforms can also be constructed on $\mathcal{L}(\Z_{2m}^{*n})$. They are $L_p$-bounded maps as well there, as shown in \cite[Theorem 3.5]{MR17}. Therefore, Theorem \ref{FreeXpChaosHilbert} can also be proved with this technique replacing $\mathbb{F}_n$ by $\Z_{2m}^{*n}$ in the statement.}  
\end{remark}

\vskip-5pt

\noindent \textbf{Acknowledgement.} We are very much indebted to \'Eric Ricard for pointing out the role of free Hilbert transforms in this context. The third-named author wants to express his gratitude to Alexandros Eskenazis and Assaf Naor for \hskip-2pt some valuable comments on a preliminary version, which led to an improved presentation. The three authors were partly supported by Grant PID 2019-107914GB-I00 (MCIN PI J. Parcet) as well as Severo Ochoa Grant CEX2019-000904-S (ICMAT), funded by MCIN/AEI 10.13039/501100011033. Jos\'e M. Conde-Alonso was also supported by the Madrid Government Program V under PRICIT Grant SI1/PJI/2019-00514. Antonio Ismael Cano-Mármol was supported by Grant SEV-2015-0554-19-3 funded by MCIN/AEI/10.13039/501100011033 and by ESF Investing in your future.

\bibliography{biblio}

\begin{thebibliography}{10}

\bibitem{Bz}
M.~Bo\.{z}ejko, T.~Januszkiewicz, and R.J. Spatzier.
\newblock Infinite {C}oxeter groups do not have {K}azhdan's property.
\newblock {\em J. Operator Theory}, 19:63--67, 1988.

\bibitem{CCP22p2}
A.I. Cano-M\'armol, J.M. Conde-Alonso, and J.~Parcet.
\newblock Trigonometric chaos and {$\mathrm{X}_p$} inequalities {II}:
  {$\mathrm{X}_p$} inequalities in group von {N}eumann algebras. {P}reprint
  2022.

\bibitem{dL}
K.~de~Leeuw.
\newblock On {$L_{p}$} multipliers.
\newblock {\em Ann. of Math.}, 81:364--379, 1965.

\bibitem{H78}
U.~Haagerup.
\newblock An example of a nonnuclear {$C^{\ast} $}-algebra, which has the
  metric approximation property.
\newblock {\em Invent. Math.}, 50:279--293, 1978/79.

\bibitem{HP93}
U.~Haagerup and G.~Pisier.
\newblock Bounded linear operators between {$C^*$}-algebras.
\newblock {\em Duke Math. J.}, 71:889--925, 1993.

\bibitem{Ho}
L.~H\"{o}rmander.
\newblock Estimates for translation invariant operators in {$L^{p}$} spaces.
\newblock {\em Acta Math.}, 104:93--140, 1960.

\bibitem{JMST79}
W.B. Johnson, B.~Maurey, G.~Schechtman, and L.~Tzafriri.
\newblock Symmetric structures in {B}anach spaces.
\newblock {\em Mem. Amer. Math. Soc.}, 19(217), 1979.

\bibitem{JM10}
M.~Junge and T.~Mei.
\newblock Noncommutative {R}iesz transforms---a probabilistic approach.
\newblock {\em Amer. J. Math.}, 132:611--680, 2010.

\bibitem{JMP14}
M.~Junge, T.~Mei, and J.~Parcet.
\newblock Smooth {F}ourier multipliers on group von {N}eumann algebras.
\newblock {\em Geom. Funct. Anal.}, 24:1913--1980, 2014.

\bibitem{JMP18}
M.~Junge, T.~Mei, and J.~Parcet.
\newblock Noncommutative {R}iesz transforms---dimension free bounds and
  {F}ourier multipliers.
\newblock {\em J. Eur. Math. Soc.}, 20:529--595, 2018.

\bibitem{JPPP17}
M.~Junge, C.~Palazuelos, J.~Parcet, and M.~Perrin.
\newblock Hypercontractivity in group von {N}eumann algebras.
\newblock {\em Mem. Amer. Math. Soc.}, 249(1183), 2017.

\bibitem{JPX07}
M.~Junge, J.~Parcet, and Q.~Xu.
\newblock Rosenthal type inequalities for free chaos.
\newblock {\em Ann. Probab.}, 35:1374--1437, 2007.

\bibitem{JX03}
M.~Junge and Q.~Xu.
\newblock Noncommutative {B}urkholder/{R}osenthal inequalities.
\newblock {\em Ann. Probab.}, 31:948--995, 2003.

\bibitem{JX08}
M.~Junge and Q.~Xu.
\newblock Noncommutative {B}urkholder/{R}osenthal inequalities. {II}.
  {A}pplications.
\newblock {\em Israel J. Math.}, 167:227--282, 2008.

\bibitem{LP86}
F.~Lust-Piquard.
\newblock In\'{e}galit\'{e}s de {K}hintchine dans {$C_p\;(1<p<\infty)$}.
\newblock {\em C. R. Acad. Sci. Paris S\'{e}r. I Math.}, 303:289--292, 1986.

\bibitem{LP98}
F.~Lust-Piquard.
\newblock Riesz transforms associated with the number operator on the {W}alsh
  system and the fermions.
\newblock {\em J. Funct. Anal.}, 155:263--285, 1998.

\bibitem{MR17}
T.~Mei and \'{E}. Ricard.
\newblock Free {H}ilbert transforms.
\newblock {\em Duke Math. J.}, 166:2153--2182, 2017.

\bibitem{Mi}
S.G. Mikhlin.
\newblock On the multipliers of {F}ourier integrals.
\newblock {\em Dokl. Akad. Nauk SSSR}.

\bibitem{N16}
A.~Naor.
\newblock Discrete {R}iesz transforms and sharp metric {$X_p$} inequalities.
\newblock {\em Ann. of Math.}, 184:991--1016, 2016.

\bibitem{NS16}
A.~Naor and G.~Schechtman.
\newblock Metric {$X_p$} inequalities.
\newblock {\em Forum Math. $\Pi$}, 4:e3, 81, 2016.

\bibitem{R70}
H.P. Rosenthal.
\newblock On the subspaces of {$L^{p}$} {$(p>2)$} spanned by sequences of
  independent random variables.
\newblock {\em Israel J. Math.}, 8:273--303, 1970.

\bibitem{S38}
I.J. Schoenberg.
\newblock Metric spaces and completely monotone functions.
\newblock {\em Ann. of Math.}, 39:811--841, 1938.

\end{thebibliography}
\bibliographystyle{plain}

\hfill \noindent \textbf{Antonio Ismael Cano M\'armol} \\
\null \hfill Instituto de Ciencias Matem{\'a}ticas 
\\ \null \hfill Consejo Superior de Investigaciones Cient{\'\i}ficas 
\\ \null \hfill\texttt{ismael.cano@icmat.es}

\vskip2pt

\hfill \noindent \textbf{Jos\'e M. Conde-Alonso} \\
\null \hfill Universidad Aut\'onoma de Madrid 
\\ \null \hfill Instituto de Ciencias Matem{\'a}ticas  
\\ \null \hfill\texttt{jose.conde@uam.es}

\vskip2pt

\hfill \noindent \textbf{Javier Parcet} \\
\null \hfill Instituto de Ciencias Matem{\'a}ticas 
\\ \null \hfill Consejo Superior de Investigaciones Cient{\'\i}ficas 
\\ \null \hfill\texttt{parcet@icmat.es}

\end{document}